\documentclass[sn-mathphys-num]{sn-jnl}
\usepackage{graphicx}
\usepackage{multirow}
\usepackage{amsmath,amssymb,amsfonts}
\usepackage{amsthm}
\usepackage{mathrsfs}
\usepackage[title]{appendix}
\usepackage{xcolor}
\usepackage{textcomp}
\usepackage{manyfoot}
\usepackage{booktabs}
\usepackage{algorithm}
\usepackage{algorithmicx}
\usepackage{algpseudocode}
\usepackage{listings}
\usepackage{subcaption}

\newtheorem{example}{Test problem}

\begin{document}

\title[{\bf Petrov-Galerkin Variational Physics-Informed Neural Network Framework for Two-Dimensional Singularly Perturbed Problems}]{\bf Petrov-Galerkin Variational Physics-Informed Neural Network Framework for Two-Dimensional Singularly Perturbed Problems}

\author[1]{\fnm{Vijay} \sur{Kumar}}\email{vijaykumarmathematics@gmail.com}

\author*[1]{\fnm{Gautam} \sur{Singh}}\email{gautam@nitt.edu}

\affil[1]{\orgdiv{Department of Mathematics}, \orgname{National Institute of Technology Tiruchirappalli}, \city{Tiruchirappalli}, \postcode{620015}, \state{Tamil Nadu}, \country{India}}

\abstract{This study proposes a Petrov-Galerkin based Variational Physics-Informed Neural Network (VPINN) for efficiently solving two-dimensional singularly perturbed problems (SPPs) with one and two small perturbation parameters. The approach employs neural networks to construct the trial solution space, while tensor-product hat functions are adopted as test functions to enforce the variational form. To accurately resolve of sharp boundary layers, the variational form is implemented using a Petrov-Galerkin formulation. Dirichlet boundary conditions are imposed directly, while the source terms are computed using automatic differentiation. Computational experiments on standard two-dimensional problems demonstrate that the proposed method achieves high accuracy in both the maximum and $L_2$  norms. These results confirm the efficiency and robustness of the Petrov-Galerkin VPINN approach in accurately capturing the multiscale features of two-dimensional SPPs.}



\maketitle

\section{Introduction}\label{secint}

SPPs arise across many branches of science and engineering, such as chemical engineering~\cite{Alhu}, control systems~\cite{kokotovic}, modeling of flow and transport phenomena~\cite{Modeling_Fluid}, water quality analysis~\cite{Water_pollution}, and fluid dynamics~\cite{Schlicting}. These problems involve small perturbation parameters that multiply the highest-order derivatives, producing boundary or interior layers where the solution exhibits rapid variation. Such features make the numerical and analytical treatment of SPPs particularly difficult.

Over the years, a wide range of numerical approaches have been developed for tackling one- and two-parameter SPPs. Among these, the finite element method (FEM) has been extensively studied and applied in various formulations~\cite{Weak_galerkin, Gautam2D,  FEM_2D}. A highly influential class within this family is the Discontinuous Galerkin (DG) framework, originally proposed by Reed and Hill for neutron transport equations~\cite{Reed}. DG methods are now widely used because they handle complex geometries effectively and provide robust PDE approximations. To efficiently manage the numerical difficulties caused by singular perturbations in two-dimensional problems, several formulations of the DG method have been developed. Among the most widely used are the Interior Penalty Discontinuous Galerkin scheme~\cite{IPDG}, the Symmetric Interior Penalty Galerkin scheme~\cite{SIPG}, and the Non-Symmetric Interior Penalty Galerkin approach~\cite{Gautam2D}. For a detailed theoretical foundation and comprehensive discussion of DG methods, one can refer ~\cite{Riviere}.

In parallel with developments in traditional numerical methods, deep neural networks (DNNs) have emerged as powerful tools for learning solutions of differential equations directly from governing physical laws. Building on this concept, Raissi et al.~\cite{PINN_Raissi} introduced the framework of Physics-Informed Neural Networks (PINNs), where the training process enforces the governing equations accompanied by the associated  initial and boundary conditions through a physics-based loss function. In contrast to traditional mesh-dependent numerical techniques, PINNs function in a mesh-free manner, offering flexibility in handling high-dimensional spaces and irregular computational domains.

To address the challenges associated with non-smooth or highly oscillatory solutions, several variational extensions of PINNs have been introduced. A Petrov-Galerkin formulation employs neural networks as nonlinear trial functions while adopting orthogonal bases, such as Legendre polynomials and trigonometric functions for the test space~\cite{VPINN1}. The hp-VPINN framework~\cite{VPINN} further advances this concept by incorporating domain decomposition and projecting onto higher-order polynomial subspaces. Within this framework, the neural network defines a global trial function, whereas the test functions are locally supported polynomials defined independently on each subdomain.

The Finite Basis PINN approach~\cite{FB_PINN} is inspired by the FEM, where the approximate solution is expressed as a weighted sum of finite basis functions. This hybrid formulation enhances the representational capacity of PINNs, enabling them to capture multi-scale structures and sharp gradients frequently encountered in SPPs. The framework was further advanced in the work of Raina et al.~\cite{FB-PINN_Natesan} which addresses the two small parameters for the convection-dominated models, achieving improved stability and accuracy compared to conventional PINN formulations. Cao et al.~\cite{PAPINN} consider the Parameter-Asymptotic PINN to improve the performance of standard PINNs in resolving steep boundary layers arising in SPPs. Subsequently, Boro et al.~\cite{PINN_Natesan} extended this concept to two small parameters in convection-dominated regimes, demonstrating that the inclusion of asymptotic information can substantially enhance the precision of the obtained solutions.

Yadav and Ganesan~\cite{Yadav2} proposed an ANN-driven stabilization mechanism for the Streamline Upwind Petrov-Galerkin FEM to improve its performance on SPPs. Additionally, DG-inspired neural network architectures that embed fundamental discontinuous Galerkin concepts have been presented in~\cite{PINN_NIPG, franck_2024, sun_DG}. 

In this study, we develop a VPINN with the Petrov-Galerkin framework to solve one- and two-parameter singularly perturbed BVPs, as well as parabolic PDEs involving small perturbation parameters. SPPs are challenging for classical numerical methods, which often require layer-adapted meshes and prior knowledge of boundary layer locations, and for standard strong-form PINNs, which suffer from instability and slow convergence in the presence of sharp gradients. The proposed approach addresses these limitations by employing a mesh-free variational formulation with neural networks as the trial space and localized tensor product of hat functions as the test space. Dirichlet boundary conditions are imposed directly, while the source terms are computed using automatic differentiation, avoiding the requirement for explicit analytical forms. For time-dependent problems, Kumar and Singh \cite{kumar2025_VPINN} adopted a hybrid strategy that combines the backward Euler scheme for temporal discretization with the VPINN formulation for SPPs in  one-spatial dimensions. However, to the best of our knowledge, a Petrov–Galerkin VPINN framework for two-dimensional SPPs has not yet been reported in the literature.

The primary contribution of this work lies in generalizing the applicability of the VPINN by addressing more complex SPPs.  The proposed algorithm can effectively handle different solution regimes without requiring prior knowledge of boundary layer locations or modifications to the neural network architecture. Furthermore, a progressive refinement strategy is employed during training to enhance the representational capability of the neural network, thereby improving both the accuracy and robustness of convergence when approximating multiscale solutions with steep gradients.

The remainder of this article is organized in the following manner: Section~\ref{sec1} introduces the fundamental concepts of ANNs, PINNs, and VPINNs, along with a brief description of the dataset used in this study. Section~\ref{sec3} presents the formulation of the model, the proposed algorithm, and representative numerical examples. Finally, the paper concludes with a summary of the key findings in Section~\ref{Conclusion}.

\section{Overview of neural network models}\label{sec1}

\subsection{Neural network preliminaries}\label{NN}
Suppose $F:\mathbb{R}^n \to \mathbb{R}^{nout}$  denote a standard feed-forward ANN consisting  $r-1$ hidden layers, where the $P$-th hidden layer contains $n_P$ neurons. The architecture mainly consists of  an input layer, multiple hidden layers, and an output layer. A commonly used configuration is the Multilayer Perceptron (MLP), a feed-forward neural network in which each neuron in a given layer is fully connected to all neurons in the subsequent layer. Depending on the depth of the architecture, an MLP with a single hidden layer is referred to as a shallow neural network, whereas models with multiple hidden layers are termed DNNs. Neural networks are broadly known to approximate nonlinear high-dimensional functions, and therefore are considered universal function approximators in scientific computing.

For a spatial input $(x, y) \in \mathbb{R}^2$, the approximation $u_\theta(x, y) \in \mathbb{R}$ is generated by a fully connected feed-forward neural network composed of successive affine transformations and nonlinear activations. For a network with $P$ layers, the forward propagation can be expressed as
\[
u_\theta(x, y) = \mathcal{N}_P \circ \Phi \circ \mathcal{N}_{P-1} \circ \cdots \circ \Phi \circ \mathcal{N}_1(x, y),
\]
where each layer transformation satisfies
\[
\mathcal{N}_l(x, y) = W_l\, \Phi(\mathcal{N}_{l-1}(x, y)) + b_l, 
\quad 1 \le l \le P, \quad 
\text{with } \mathcal{N}_0(x, y) = (x, y).
\]
Here, $W_l \in \mathbb{R}^{N_l \times N_{l-1}}$ and $b_l \in \mathbb{R}^{N_l}$ denote the weight matrix and bias vector of the $l$-th layer. Typical activation choices for $\Phi$ include $\tanh$, sigmoid, SiLU, ReLU, and Leaky ReLU. The entire set of trainable parameters is compactly written as $\theta = \{\, W_l,\, b_l \,\}_{l=1}^{P}.$

In conventional supervised learning, neural network parameters are typically optimized using minimizing a loss function defined on labeled datasets. Whereas, scientific machine learning applications often face a scarcity of labeled data, requiring the incorporation of physical knowledge or governing equations to guide the learning process. In such cases, prior knowledge of governing physical laws can be embedded directly into the loss functional, ensuring that the neural network solution adheres to the underlying equations. 

\begin{figure}[htbp]
	\centering
	\includegraphics[width=0.7\textwidth]{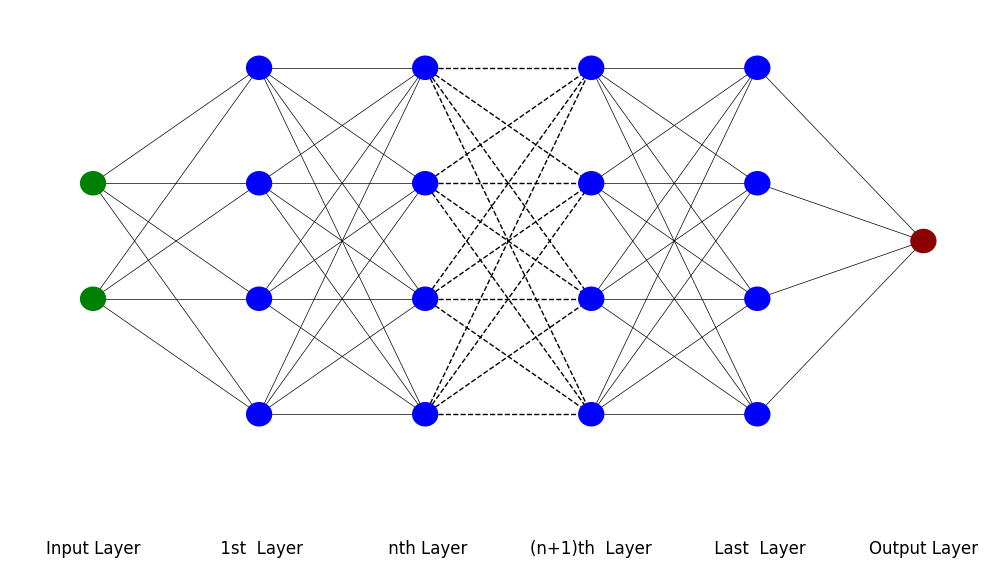}
	\label{fig:ANNs}
	\caption{Neural network visualization.}\label{figureANNs}
\end{figure}

ANNs have attracted considerable interest in computational science because of their capability to approximate complex functions universally and the adaptability provided by hyperparameter tuning. In the context of two-dimensional PDEs, ANNs have shown strong potential for capturing complex solution behaviors, particularly when conventional numerical techniques face challenges in resolving multiscale features. Several works have explored the application of ANNs for solving PDEs~\cite{DNN_FEM_1, DNN_FEM_2}, as well as hybrid approaches that combine neural networks with classical discretization methods~\cite{DNN_FEM_3}. For instance, the study in~\cite{DNN_FEM_4} proposed an enhanced FEM framework that leverages machine learning to approximate element-level mappings, thereby reducing the computational cost of large-scale simulations. Nevertheless, in such hybrid methods, the FEM is not directly embedded into the neural network structure, but instead serves as a tool for data generation or acceleration of specific computations.  Figure~\ref{figureANNs} illustrates the architecture of the two-dimensional ANNs.

\subsection{Physics-informed neural networks for 2D problems}\label{PINN}

PINNs~\cite{PINN_Raissi} are neural networks that incorporate physical laws expressed by PDEs into their training. They can be used both for data-driven solution of PDEs and for the discovery of governing equations, and they naturally extend to two-dimensional problems. Their main strength lies in the ability to embed physical laws expressed as PDEs, along with boundary and initial conditions within the training process of a neural network. This approach reduces the reliance on dense measurement data or fine spatial discretization, making PINNs particularly appealing for problems where high-fidelity data acquisition is costly or impractical. The mesh-free formulation allows the PDE solution task to be reformulated as an optimization problem, where the loss function integrates both data fidelity  and physics consistency. Recent developments include the integration of established physical models, such as rolling force formulations, as additional constraints within the loss function or the network architecture~\cite{luo2025}.

For a general two-dimensional PDE defined on a domain \( \Omega \subset \mathbb{R}^2 \), the problem may be expressed as:
\[
\mathcal{D}(u(x,y)) = r(x,y), \quad (x,y) \in \Omega, 
\qquad \mathcal{B}(u(x,y)) = 0, \quad (x,y) \in \partial \Omega,
\]
where \( \mathcal{D} \) and \( \mathcal{B} \) represent the differential operator and the boundary conditions with the given source term \( r(x,y) \). A neural network \( u_\theta(x,y) \) is characterized by a set of trainable parameters \( \theta \), which is employed to approximate the unknown solution of the governing PDEs. 
The residual of the differential operator is expressed as:
\[
v(x,y) = \mathcal{D}(u_\theta(x,y)) - r(x,y).
\]
To train the network, a composite loss functional is constructed that enforces both the physical constraints in the interior of the domain and the prescribed boundary data. 
In two spatial dimensions, the loss function is formulated as:
\[
\text{Loss}_{\text{PINN}}(\theta) 
= \frac{1}{N_f} \sum_{i=1}^{N_f} \left| v(x_i^f,y_i^f) \right|^2 
+ \frac{1}{N_u} \sum_{i=1}^{N_u} \left| u_\theta(x_i^u,y_i^u) - g(x_i^u,y_i^u) \right|^2,
\]
where the set $\{(x_i^f, y_i^f)\}_{i=1}^{N_f} \subset \Omega$ denotes the collocation points inside the computational domain where the PDE residual is minimized, and $\{(x_i^u, y_i^u)\}_{i=1}^{N_u} \subset \partial\Omega$ represents the boundary nodes at which the Dirichlet boundary condition $u(x,y) = g(x,y)$ is enforced.

As in the one-dimensional case, hard imposition of boundary conditions can be achieved by designing the trial solution in a way that satisfies them exactly, i.e., constructing the neural network architecture such that $u_\theta(x,y)\big|_{\partial \Omega} = g(x,y).$ This approach ensures the boundary conditions are inherently satisfied within the solution space, improving numerical stability and accuracy.

Despite these advantages, solving two-dimensional SPPs with PINNs remains challenging. The solution may exhibit sharp gradients or thin layers near boundaries or within the interior, which are difficult to capture for a single global neural network. To mitigate this, domain decomposition and multi-network strategies have been investigated~\cite{Deep_decomposition}, though they often depend on a priori knowledge of where layers or singularities occur. Developing robust PINN formulations for 2D SPPs without such prior information remains an active area of research.

\subsection{Variational physics-informed neural networks for 2D problems}\label{VPINN}

In two spatial dimensions, VPINNs generalize the PINN framework by enforcing the residual of the governing PDE in a  weak sense. Instead of minimizing the residual 
$\mathcal{D}(u_\theta(x,y)) - f(x,y)$ pointwise, VPINNs impose orthogonality with respect to a family of test functions $\{v_{ij}(x,y)\}_{i,j=1}^M$, i.e.,
\[
\int_{\Omega} \big( \mathcal{D}(u_\theta(x,y)) - f(x,y) \big) v_{ij}(x,y) = 0 
\quad \forall i,j = 1, \ldots, M,
\]
where $\Omega \subset \mathbb{R}^2$ denotes the computational domain.  

The corresponding VPINN loss functional is formulated as:
\[
\text{Loss}_{\text{VPINN}}(\theta) 
= \frac{1}{M^{2}}\sum_{i,j=1}^{M} \left| \int_{\Omega} 
\big( \mathcal{D}(u_\theta(x,y)) - f(x,y) \big) v_{ij}(x,y)  \right|^2 
+ \text{MSE}_u,
\]
with the boundary loss term:
\[
\text{MSE}_u = \frac{1}{N_u} \sum_{i=1}^{N_u} 
\left| u_\theta(x_i^u,y_i^u) - g(x_i^u,y_i^u) \right|^2.
\]
Alternatively, boundary conditions may be imposed as hard constraints by constructing the neural network such that $u_\theta(x,y) = g(x,y)$ holds exactly on $\partial \Omega$.  

In our implementation, we adopt tensor-product hat functions as the test space. A tensor-product hat function in 2D is formed by multiplying two 1D hat functions,
\[
v_{ij}(x,y) = \psi_i(x) \, \psi_j(y),
\]
where $\psi_i$ and $\psi_j$ are piecewise-linear basis functions associated with grid nodes in the $x$- and $y$-directions. This construction yields a locally supported bilinear basis, similar to a standard finite element test space. All integrals in the weak formulation are approximated numerically using collocation points serving as integration nodes. 

\subsection{Dataset for steady and parabolic 2D problems }\label{data}

This subsection examines datasets generated from two-dimensional SPPs. Within the VPINN framework, we employ $10 \times 10$ tensor-product hat test functions in conjunction with $40 \times 40$ quadrature nodes to accurately perform the weak-form integration. Dirichlet boundary conditions are embedded directly into the network output to ensure exact satisfaction of the prescribed values. The neural network architecture comprises $4$ hidden layers, each containing $20$ neurons, with the $\tanh$ function serving as the activation function. For both single-parameter and double-parameter two-dimensional SPPs, the L-BFGS optimizer is employed to train the model for $1000$ epochs.

For problems with temporal evolution, the same spatial resolution of $10 \times 10$ test functions is maintained, while the quadrature rule is extended into the space–time domain with $40 \times 40 \times 40$ integration nodes. The training process begins with the Adam optimizer (learning rate $10^{-3}$) for stable initialization, followed by refinement using the L-BFGS method for an additional 1000 epochs. The overall network architecture and dataset structure remain consistent with the steady-state formulation. All simulations were implemented in Python using TensorFlow and executed on a Dell workstation with an Intel\textsuperscript{\textregistered} Core\textsuperscript{\texttrademark} i$7-7700$ CPU ($3.60~$GHz), a $64-$bit operating system, and $32.0~$GB RAM.

\section{Test problems and computational experiments}\label{sec3}

\subsection{Two-dimensional one-parameter elliptic boundary value problems (BVPs)}\label{sec3.1}
Consider the two-dimensional one-parameter elliptic BVPs:
\begin{equation}\label{E1}
	\left\{
	\begin{array} {ll}
		-\epsilon \Delta u + \boldsymbol{\beta}(x,y)\cdot\nabla u + \sigma(x,y) u = f, & (x,y) \in \Omega=(0,\,1)\times(0,\,1), \\[8pt]
		u = 0, &  (x,y) \in \partial\Omega.\\[8pt]
	\end{array}\right.
\end{equation}
Here, $0<\epsilon\ll 1$ is a small parameter, $\boldsymbol{\beta}(x,y) = (\beta_1, \beta_2)^{T}$ with $\beta_1 \geq  \alpha_1 > 0$, $\beta_2 \geq  \alpha_2 > 0$,
and $\sigma(x,y) \geq \delta > 0$ for any $(x,y)\in\overline{\Omega}$. 
Furthermore, we assume the stability condition:
\[c_0^2 = \left(\sigma - \tfrac{1}{2}\nabla\cdot \boldsymbol{\beta}\right) \geq \omega > 0.\]

We assume that $ \boldsymbol{\beta}, \sigma$ and $f$ are smooth enough and $f$ meets the following compatibility property:
\[f(0,0)=f(0,1)=f(1,0)=f(1,1)=0.\]

The weak form corresponding to the BVP~\eqref{E1} is to determine a function 
\( u \in H_0^1(\Omega) \) satisfying, for every test function \( v \in H_0^1(\Omega) \),
\begin{equation}\label{E2}
	B(u,v) = (f,v),
\end{equation}
here the bilinear operator \( B(\cdot,\cdot) \) is defined as:
\[
B(u,v) 
= \int_{\Omega} \epsilon \, \nabla u \!\cdot\! \nabla v \, \mathrm{d}\Omega
+ \int_{\Omega} \big( \boldsymbol{\beta}\!\cdot\!\nabla u + \sigma u \big) v \, \mathrm{d}\Omega,
\]
the right-hand side of the variational formulation corresponds to the linear functional
\[
(f,v) = \int_{\Omega} f\,v \, \mathrm{d}\Omega.
\]

Let \( \Omega_h \) denote a conforming partition of the physical domain \( \Omega \subset \mathbb{R}^2 \) into a finite collection of elements. Based on this discretization, the associated finite-dimensional subspace is \( V_h \subset H_0^1(\Omega) \), spanned by appropriate basis functions. 

In the neural formulation, the trial space is represented by a parameterized DNN,
\[
U_{\text{NN}} 
= \left\{ u_\theta(x,y) = \mathrm{DNN}(x,y; \theta) 
\,\middle|\, \theta \in \mathbf{\omega} \right\},
\]
where \( \mathrm{DNN}(x,y;\theta) \) represents a composition of affine transformations and nonlinear activation function,  and \( \mathbf{\omega} \) denotes the set of trainable parameters defining the network.

The discrete variational formulation of problem~\eqref{E2} can now be stated as follows:
\[
\text{Find } u_h \in U_{NN} \text{ such that } 
B_h(u_h, v) = (f, v), \quad \forall v \in V_h.
\]
Here, \( B_h(\cdot, \cdot) \) denotes the bilinear form evaluated elementwise over the mesh \( \Omega_h \).

\begin{figure}[htbp]
	\centering
	\includegraphics[width=1\textwidth]{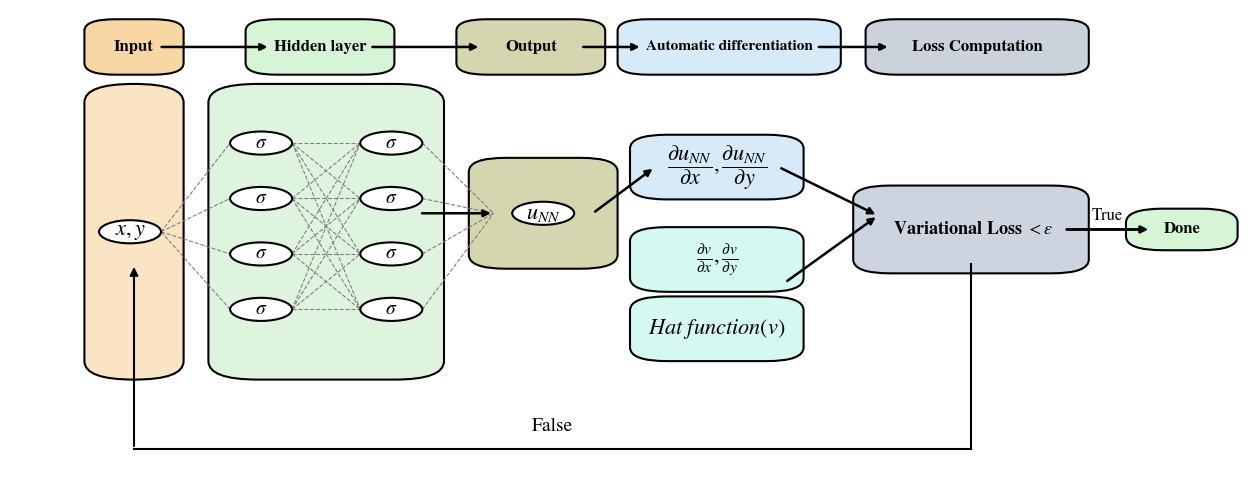}
	\caption{VPINN algorithm flowchart.}\label{figure_algorithm1}
\end{figure}

The VPINN procedure for solving two-dimensional one-parameter elliptic BVPs is summarized in Algorithm~\ref{Algorithm1}, and its graphical representation is illustrated in Figure~\ref{figure_algorithm1}.

\begin{algorithm}
	\caption{VPINN framework for solving $-\epsilon \Delta u + \beta_1(x,y)u_x + \beta_2(x,y)u_y + \sigma(x,y)u = f(x,y)$ in $\Omega=(0,1)^2$ with homogeneous Dirichlet boundary conditions.}
	\label{alg:vpinn-2d-lbfgs}
	\begin{algorithmic}[1]
		\State \textbf{Input:} Diffusion coefficient $\epsilon > 0$, consider $N$ collocation points and $M$ test functions in the computational domain.
		\State \textbf{Construct trial solution:} $u_\theta(x,y) = x(1-x)y(1-y)\,\text{NN}_\theta(x,y)$, ensuring $u_\theta=0$ on $\partial\Omega$.
		\State \textbf{Generate:} Uniformly distributed collocation points $\{(x_i,y_j)\}_{i,j=1}^{N}$ within the interior of $\Omega$.
		\State \textbf{Formulate test space:} Define a family of two-dimensional basis functions $\{v_{kl}(x,y)\}_{k,l=1}^{M}$ by tensorizing 1D piecewise linear hat functions.
		\State \textbf{Evaluate forcing term:} Compute $f(x,y) = -\epsilon(u_{xx}+u_{yy}) + \beta_1 u_x + \beta_2 u_y + \sigma u$ using automatic differentiation of the exact solution $u_{\text{exact}}$.
		\For{each $v_{kl}$ in the test space}
		\State \textbf{Assemble weak residual:}
		\[
		R_{kl}(\theta) = \int_{\Omega}\!\!\big[-\epsilon(\nabla^2 u_\theta) + \beta_1 u_{\theta,x} + \beta_2 u_{\theta,y} + \sigma u_\theta - f\big] v_{kl}(x,y)\,d\Omega.
		\]
		\EndFor
		\State \textbf{Loss functional formulation:}
		\[
		\mathcal{L}(\theta) = \frac{1}{M^2}\sum_{k,l=1}^{M}\big(R_{kl}(\theta)\big)^2.
		\]
		\State \textbf{Initialize:} Neural network parameters $\theta$.
		\State \textbf{Optimization:} Minimize $\mathcal{L}(\theta)$ via the L-BFGS algorithm with gradient updates computed through automatic differentiation.
		\State \textbf{Return:} Trained neural approximation $u_\theta(x,y)$ and computed error norms $\|u_\theta - u_{\text{exact}}\|_{L^2}$ and $\|u_\theta - u_{\text{exact}}\|_{L^\infty}$.
	\end{algorithmic}\label{Algorithm1}
\end{algorithm}

Two test problems are considered: one involving constant coefficients and another involving variable coefficients.

\begin{example}\label{example1}
	Consider the $2$D one-parameter elliptic BVP:
	\begin{equation}\label{E3}
		\left\{
		\begin{array} {ll}
			-\epsilon \Delta u - u_x- u_y +  u = f, & (x,y) \in \Omega= (0, 1)\times (0,1), \\[8pt]
			u = 0, &  (x,y) \in \partial\Omega.\\[8pt]
		\end{array}\right.
	\end{equation}
	The source function $f$ is derived from the exact solution
	\begin{eqnarray}
		u(x,y)=\sin(1-x)\sin(1-y)(1-e^{-x/\epsilon})(1-e^{-y/\epsilon}).\nonumber
	\end{eqnarray}
\end{example}

\begin{example}\label{example2}
	Consider the $2$D one-parameter elliptic BVP:
	\begin{equation}\label{E4}
		\left\{
		\begin{array} {ll}
			-\epsilon \Delta u +(1-x) u_x +(1-y) u_y +  u = f, & (x,y) \in \Omega= (0, 1)\times (0,1), \\[8pt]
			u = 0, &  (x,y) \in \partial\Omega.\\[8pt]
		\end{array}\right.
	\end{equation}
	The source function $f$ is derived from the exact solution
	\begin{eqnarray}
		u(x,y)=\sin(x)\sin(y)(1-e^{-(1-x)/\epsilon})(1-e^{-(1-y)/\epsilon}).\nonumber
	\end{eqnarray}
\end{example}

The algorithm is executed for $1000$ iterations. The resulting maximum error, $L_2$ error, and loss at the final iteration are reported in Table~\ref{tbl_1} for test problem~\ref{example1} and Table~\ref{tbl_2} for test problem~\ref{example2}.

\begin{figure}[htbp]
	\centering
	\begin{subfigure}[b]{0.9\textwidth}
		\includegraphics[width=\textwidth]{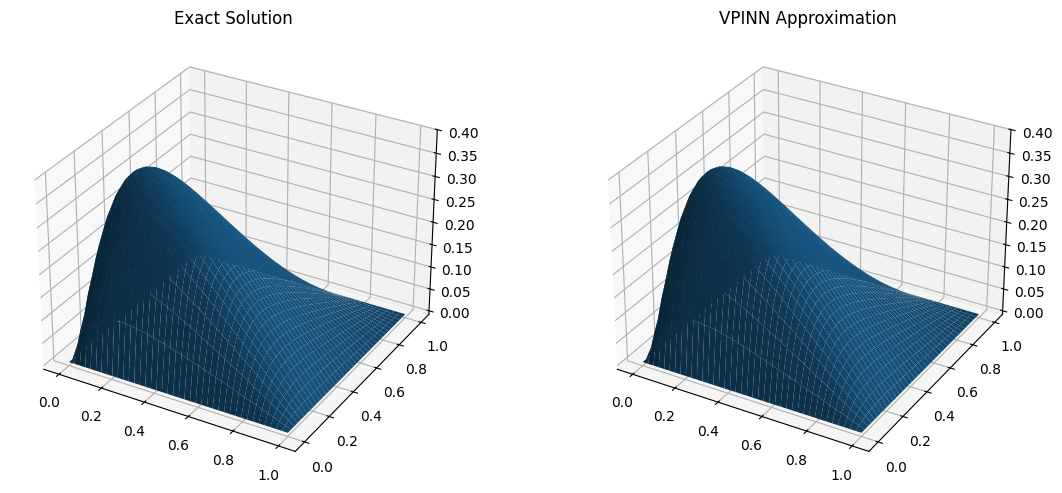}
		\caption{ $\epsilon=10^{-1}$}
		\label{fig:subfig1A}
	\end{subfigure}
	\hfill
	\begin{subfigure}[b]{0.9\textwidth}
		\includegraphics[width=\textwidth]{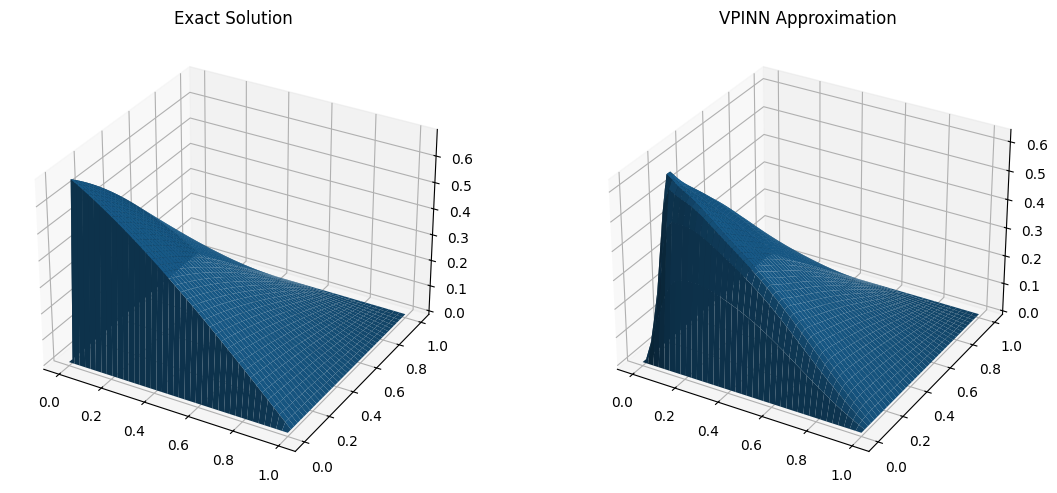}
		\caption{  $\epsilon=10^{-3}$}
		\label{fig:subfig1B}
	\end{subfigure}
	\hfill
	\begin{subfigure}[b]{0.9\textwidth}
		\includegraphics[width=\textwidth]{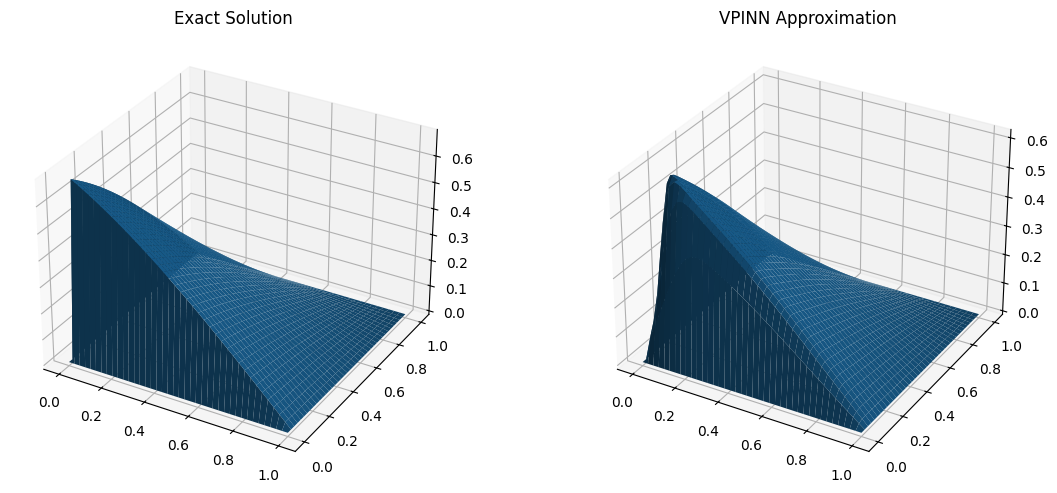}
		\caption{  $\epsilon=10^{-5}$}
		\label{fig:subfig1C}
	\end{subfigure}
	\caption{Comparison of the analytical and neural network solutions for test problem~\ref{example1} corresponding to various values of the perturbation parameters.}\label{figure1}
\end{figure}

The comparison between the  analytical and neural network solutions for different parameter values is shown in Figure~\ref{figure1} for test problem~\ref{example1} and Figure~\ref{figure4} for test problem~\ref{example2}.  

\begin{figure}[htbp]
	\centering
	
	\begin{subfigure}[b]{0.45\textwidth}
		\includegraphics[width=\textwidth,]{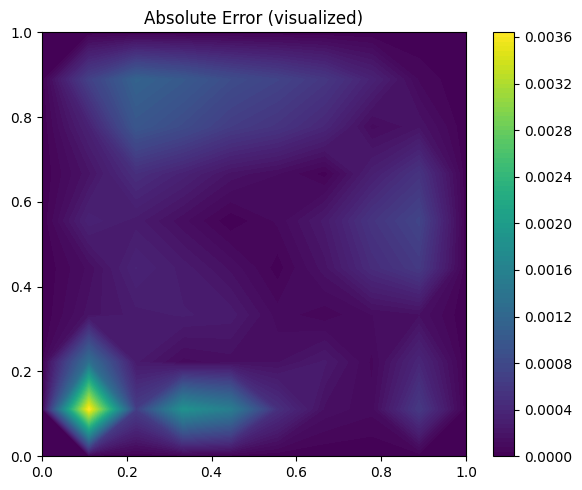}
		\caption{  $\epsilon=10^{-1}$}
		\label{fig:subfig2A}
	\end{subfigure}
	\hfill
	\begin{subfigure}[b]{0.45\textwidth}
		\includegraphics[width=\textwidth,]{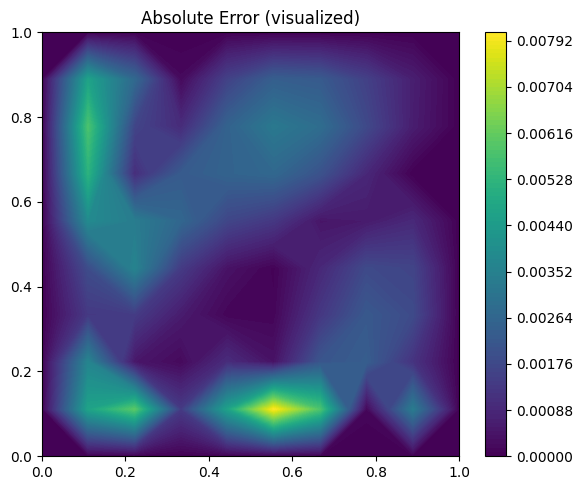}
		\caption{  $\epsilon=10^{-3}$}
		\label{fig:subfig2B}
	\end{subfigure}
	\hfill
	\begin{subfigure}[b]{0.45\textwidth}
		\includegraphics[width=\textwidth,]{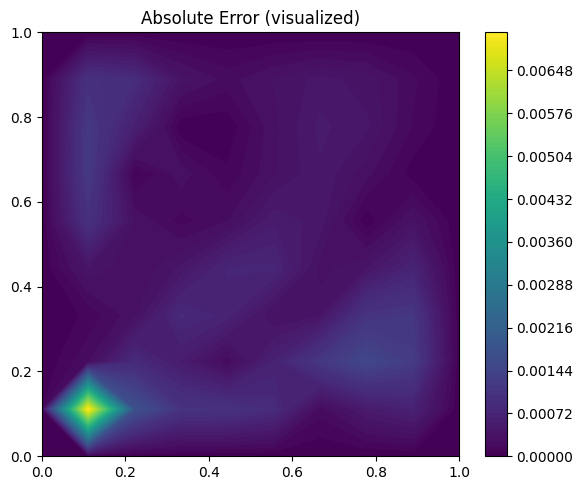}
		\caption{  $\epsilon=10^{-5}$}
		\label{fig:subfig2C}
	\end{subfigure}
	\caption{Absolute error for test problem~\ref{example1} corresponding to various values of the perturbation parameters.}\label{figure2}
\end{figure}

\begin{figure}[htbp]
	\centering
	
	\begin{subfigure}[b]{0.32\textwidth}
		\includegraphics[width=\textwidth,]{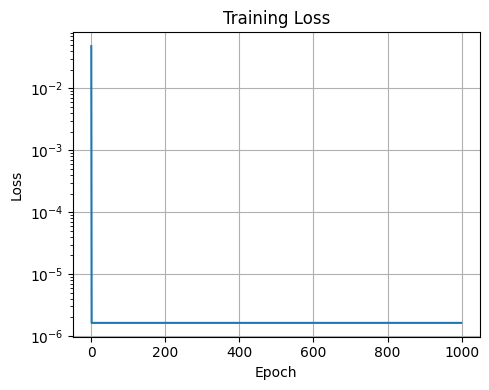}
		\caption{$\epsilon=10^{-1}$}
		\label{fig:subfig3A}
	\end{subfigure}
	\hfill
	\begin{subfigure}[b]{0.32\textwidth}
		\includegraphics[width=\textwidth,]{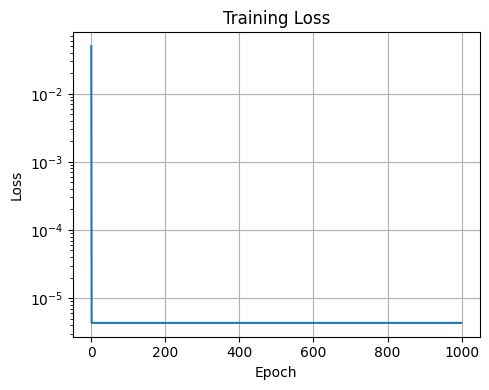}
		\caption{$\epsilon=10^{-3}$}
		\label{fig:subfig3B}
	\end{subfigure}
	\hfill
	\begin{subfigure}[b]{0.32\textwidth}
		\includegraphics[width=\textwidth,]{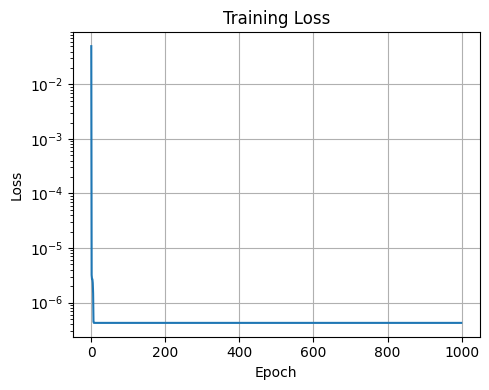}
		\caption{$\epsilon=10^{-5}$}
		\label{fig:subfig3C}
	\end{subfigure}
	\caption{Training loss versus epoch for test problem~\ref{example1} corresponding to  various values of the perturbation parameters.}\label{figure3}
\end{figure}

\begin{figure}[htbp]
	\centering
	\begin{subfigure}[b]{0.9\textwidth}
		\includegraphics[width=\textwidth]{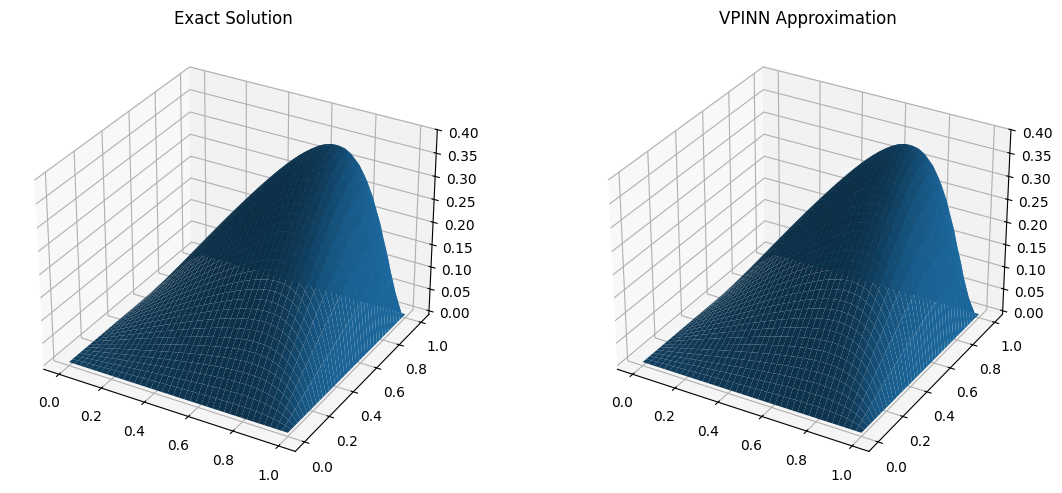}
		\caption{  $\epsilon=10^{-1}$}
		\label{fig:subfig4A}
	\end{subfigure}
	\hfill
	\begin{subfigure}[b]{0.9\textwidth}
		\includegraphics[width=\textwidth]{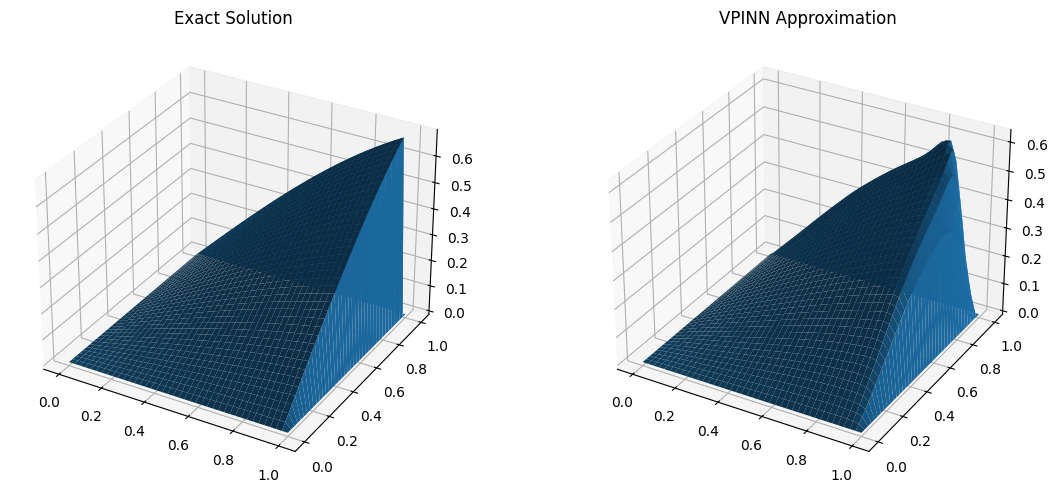}
		\caption{  $\epsilon=10^{-3}$}
		\label{fig:subfig4B}
	\end{subfigure}
	\hfill
	\begin{subfigure}[b]{0.9\textwidth}
		\includegraphics[width=\textwidth]{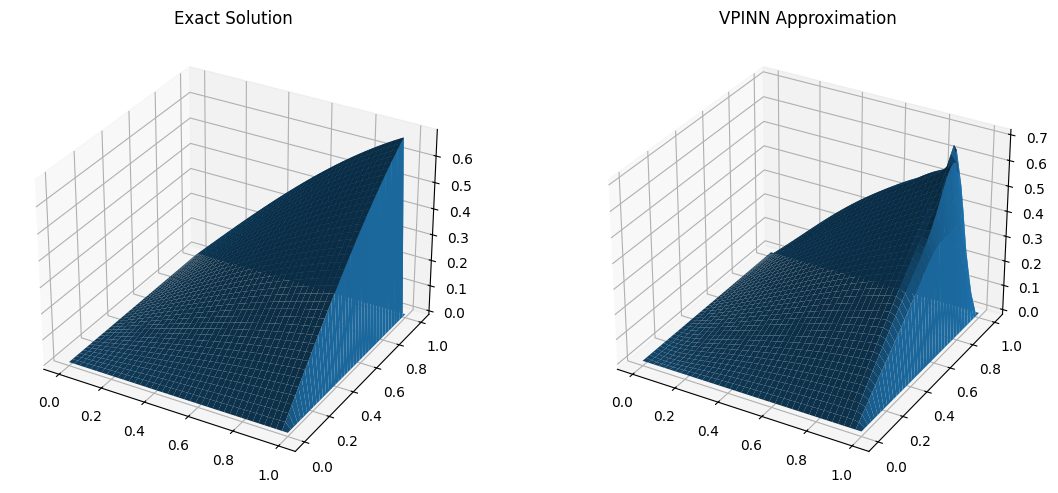}
		\caption{  $\epsilon=10^{-5}$}
		\label{fig:subfig4C}
	\end{subfigure}
	\caption{Comparison of the analytical and neural network solutions for test problem~\ref{example2} corresponding to various values of the perturbation parameters.}\label{figure4}
\end{figure}

\begin{figure}[htbp]
	\centering
	
	\begin{subfigure}[b]{0.45\textwidth}
		\includegraphics[width=\textwidth,]{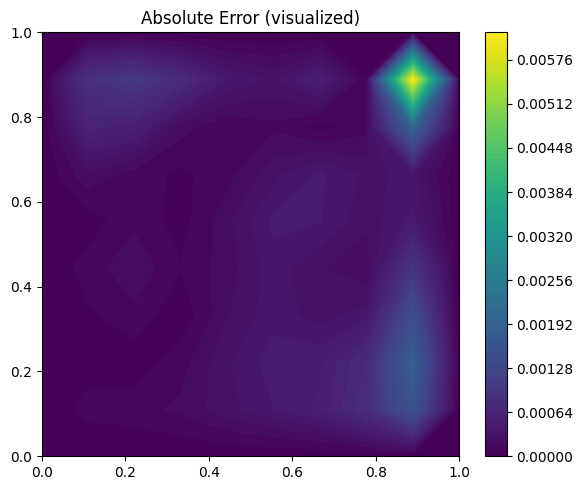}
		\caption{  $\epsilon=10^{-1}$}
		\label{fig:subfig5A}
	\end{subfigure}
	\hfill
	\begin{subfigure}[b]{0.45\textwidth}
		\includegraphics[width=\textwidth,]{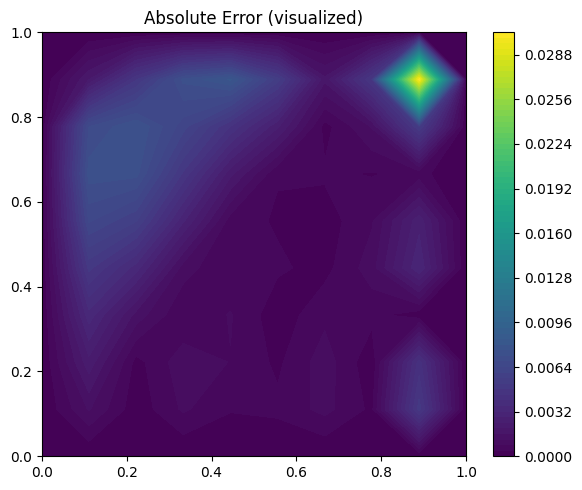}
		\caption{  $\epsilon=10^{-3}$}
		\label{fig:subfig5B}
	\end{subfigure}
	\hfill
	\begin{subfigure}[b]{0.45\textwidth}
		\includegraphics[width=\textwidth,]{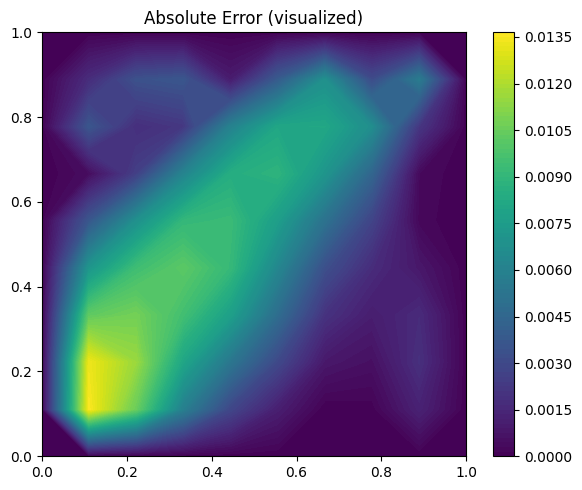}
		\caption{  $\epsilon=10^{-5}$}
		\label{fig:subfig5C}
	\end{subfigure}
	\caption{Absolute error for test problem~\ref{example2} corresponding to various values of the perturbation parameters.}\label{figure5}
\end{figure}

\begin{figure}[htbp]
	\centering
	
	\begin{subfigure}[b]{0.32\textwidth}
		\includegraphics[width=\textwidth,]{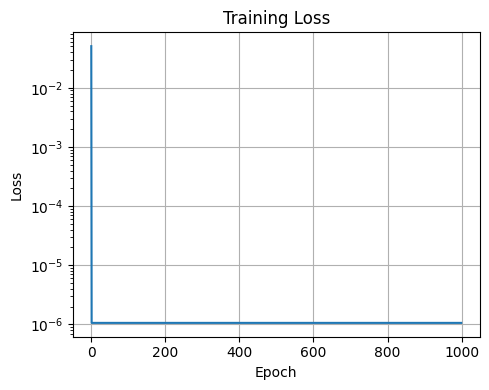}
		\caption{$\epsilon=10^{-1}$}
		\label{fig:subfig6A}
	\end{subfigure}
	\hfill
	\begin{subfigure}[b]{0.32\textwidth}
		\includegraphics[width=\textwidth,]{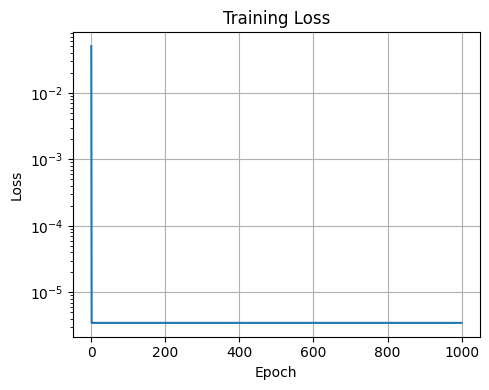}
		\caption{$\epsilon=10^{-3}$}
		\label{fig:subfig6B}
	\end{subfigure}
	\hfill
	\begin{subfigure}[b]{0.32\textwidth}
		\includegraphics[width=\textwidth,]{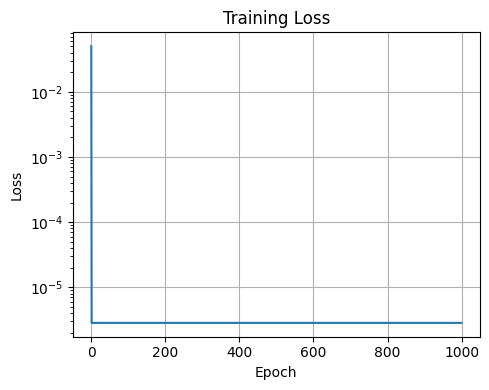}
		\caption{$\epsilon=10^{-5}$}
		\label{fig:subfig6C}
	\end{subfigure}
	\caption{Training loss versus epoch for test problem~\ref{example2} corresponding to various values of the perturbation parameters.}\label{figure6}
\end{figure}

The absolute error plots are presented in Figure~\ref{figure2} for test problem~\ref{example1} and Figure~\ref{figure5} for test problem~\ref{example2}, while the  training loss versus epoch curves corresponding to the final iteration are shown in Figure~\ref{figure3} for test problem~\ref{example1} and Figure~\ref{figure6} for test problem~\ref{example2}.

\begin{table}[htbp]
	\caption{\it{ Illustration of training loss and numerical error evolution for test problem~\ref{example1}.}}\label{tbl_1}  
	\begin{tabular}{@{}lllllll@{}}
		\multicolumn 1 {c}{}  & \multicolumn 5 {c}
		{}\\
		\hline
		$\epsilon$ &Training Loss & Max Error & Rel-Max Error & $L_2$ Error & Rel-$L_2$ Error\\
		\hline
		$10^{-1}$ & 1.6221e-06  & 3.6348e-03  & 9.2855e-03 & 5.6646e-04 & 3.4520e-03 \\ [6pt]
		
		$10^{-2}$  & 1.3240e-06 & 3.7292e-03  & 6.1872e-03 & 1.0763e-03 & 5.1047e-03 \\ [6pt]
		
		$10^{-3}$   &  4.3227e-06 & 8.0425e-03   & 1.3343e-02 & 2.1284e-03 & 1.0095e-02  \\ [6pt]
		
		$10^{-4}$  & 1.2243e-06 & 7.3734e-03  &  1.2233e-02 & 2.1658e-03  & 1.0272e-02 \\ [6pt]

		$10^{-5}$  &  4.2315e-07   &  7.0794e-03  & 1.1745e-02 & 8.8072e-04 & 4.1771e-03  \\ [6pt]
		
		\hline
	\end{tabular}
\end{table}

\begin{table}[htbp]
	\caption{\it{ Illustration of training loss and numerical error evolution for test problem~\ref{example2}.}}\label{tbl_2}  
	\begin{tabular}{@{}lllllll@{}}
		\multicolumn 1 {c}{}  & \multicolumn 5 {c}
		{}\\
		\hline
		$\epsilon$ &Training Loss & Max Error & Rel-Max Error & $L_2$ Error & Rel-$L_2$ Error\\
		\hline
		$10^{-1}$ & 1.0644e-06 & 6.1038e-03 & 1.5593e-02 & 7.6003e-04 & 4.6316e-03 \\ [6pt]
		
		$10^{-2}$  &  9.0076e-07 & 4.6403e-03  & 7.6988e-03 &  1.0683e-03 &  5.0667e-03 \\ [6pt]
		
		$10^{-3}$   &  3.4797e-06 & 3.0321e-02   & 5.0304e-02 & 4.0710e-03 & 1.9308e-02   \\ [6pt]
		
		$10^{-4}$  & 2.9361e-06 &  3.4386e-02  &  5.7048e-02 & 3.9989e-03  &  1.8966e-02 \\ [6pt]

		$10^{-5}$  & 2.8475e-06   & 1.3575e-02  & 2.2521e-02 & 4.9764e-03 & 2.3602e-02  \\ [6pt]
		
		\hline
	\end{tabular}
\end{table}

\subsection{Two-dimensional parabolic one-parameter initial boundary-value problems (IBVPs) }\label{sec3.2}
Consider the  two-dimensional one-parameter singularly perturbed parabolic IBVPs defined on the domain $G$ = $\Omega$ $\times$ $(0, T]$, with $\Omega = (0, 1)\times (0,1):$
\begin{equation}\label{E5}
	\left\{
	\begin{array} {ll}
		\frac{\partial u}{\partial t}-\epsilon \Delta u + \boldsymbol{\beta}(x,y)\cdot\nabla u + \sigma(x,y) u = f(x,y,t), & (x,y,t) \in G, \\[8pt]
		u(x,y,0)=u_0(x,y), &  (x,y) \in \bar{\Omega},\\[8pt]
		u(x,y,t)=0 & (x,y,t) \in \partial \Omega\times(0,T].
	\end{array}\right.
\end{equation}

Here, $0 < \epsilon \ll 1$ denotes a small perturbation parameter. The advection vector $\boldsymbol{\beta}(x,y) = [\beta_1(x,y), \beta_2(x,y)]^{T}$ satisfies $\beta_1 \ge \alpha_1 > 0$ and $\beta_2 \ge \alpha_2 > 0$, while the reaction coefficient satisfies $\sigma(x,y) \ge \delta > 0$ for all $(x,y) \in \overline{\Omega}$. 
Furthermore, we assume the stability condition:
\[
c_0^2 = \left(\sigma - \tfrac{1}{2}\nabla\cdot \boldsymbol{\beta}\right) \geq \omega > 0.
\]

The weak form corresponding to the BVP~\eqref{E5} is to determine a function 
\( u \in H_0^1(\Omega) \) satisfying, for every test function \( v \in H_0^1(\Omega) \),
\begin{equation}\label{E6}
	B(u,v) = (f,v),
\end{equation}
where the bilinear operator \( B(\cdot,\cdot) \) is defined as:
\[ \int_{\Omega}\frac{\partial{u}}{\partial{t}}v+\int_{\Omega} \epsilon\, \nabla u \cdot \, \nabla v  
+ \int_{\Omega} (\boldsymbol{\beta}\cdot\nabla u + \sigma u) v= \int_{\Omega}f(x,y,t)v. \]
At each discrete time level $t = t_m$, where $1 \leq m \leq M$, we obtain
\begin{equation}\label{E7}
	B(u^{m},v)=(f^{m},v),
\end{equation}
\[ \int_{\Omega}\frac{\partial{u^{m}}}{\partial{t}}v+\int_{\Omega} \epsilon\, \nabla u^{m} \cdot \, \nabla v  
+ \int_{\Omega} (\boldsymbol{\beta}\cdot\nabla u^{m} + \sigma u^{m}) v= \int_{\Omega}f^{m}v \]
\[ \int_{\Omega}\frac{(u^{m}-u^{m-1})}{\Delta t}v+\int_{\Omega} \epsilon\, \nabla u^{m} \cdot \, \nabla v  
+ \int_{\Omega} (\boldsymbol{\beta}\cdot\nabla u^{m} + \sigma u^{m}) v= \int_{\Omega}f^{m}v \]
\[C(u^{m},v)=(r^{m},v)\]
\[C(u^{m},v)=\int_{\Omega} \epsilon\, \nabla u^{m} \cdot \, \nabla v  
+ \int_{\Omega} \boldsymbol{\beta}\cdot\nabla u^{m}v +\bigg(\sigma+ \frac{1}{\Delta t}\bigg)u^{m}v\]
\[(r^{m},v)=\int_{\Omega}\bigg(f^{m}+\frac{u^{m-1}}{\Delta t}\bigg)v.\]

We introduce a finite element mesh $\Omega_h$ that partitions the computational domain $\Omega$ in a finite number of non-overlapping elements. The corresponding discrete variational formulation of equation~\eqref{E7} is given by:

\begin{equation*}
	\text{Find } u_h \in U_{NN} \text{ such that } 
	C_h(u_h^{m}, v) = (r^{m}, v), \quad \forall v \in V_h,
\end{equation*}
where $C_h(\cdot, \cdot)$ denotes element-wise integration over the mesh $\Omega_h$. The complete VPINN procedure for solving the two-dimensional parabolic PDEs is summarized in Algorithm~\ref{Algorithm2}.

\begin{algorithm}[H]
	\caption{VPINN framework for solving $u_t - \epsilon \Delta u + \beta_1(x,y) u_x + \beta_2(x,y) u_y + \sigma(x,y) u= f(x,y,t)$ on $(x,y,t) \in \Omega \times (0,1]$, where $\Omega=(0,1)\times(0,1)$, with $u=0$ on $\partial \Omega\times(0,T]$.}
	\label{Algorithm2}
	\begin{algorithmic}[1]
		\State \textbf{Input:} Diffusion parameter $\epsilon > 0$,  $\boldsymbol{\beta}=(\beta_1,\beta_2)$, reaction $\sigma$, final time $T$, spatial points $N_x,N_y$, time steps $N_t$, $M$ test functions are considered.
		\State \textbf{Construct trial function:} $u_\theta(x,y) := x(1-x)\,y(1-y)\, \text{NN}_\theta(x,y)$ to enforce Dirichlet BCs.
		\State \textbf{Generate:} Uniformly distributed collocation points $\{(x_i,y_j)\}_{i,j=1}^{N_x,N_y} \subset (0,1)\times(0,1)$, with time grid $t_n = n \Delta t,\;\Delta t = \tfrac{1}{N_t}.$
		\State \textbf{Construct:} Tensor-product hat test functions $\{v_{kl}(x,y)\}_{k,l=1}^M$ on a uniform mesh.
		\State \textbf{Compute:} Manufactured source term 
		\[
		f(x,y,t) = u_t - \epsilon(u_{xx}+u_{yy}) + \beta_1(x,y) u_x + \beta_2(x,y) u_y + \sigma(x,y) u
		\] 
		using autograd from exact solution $u_\text{exact}(x,y,t).$
		\State \textbf{Initialize:} $u^0(x,y) = u_\text{exact}(x,y,0)$.
		\For{$n=1$ to $N_t$}
		\State \textbf{Form residual at time $t_n$:}
		\[
		R(x,y) = \frac{u_\theta(x,y) - u^{n-1}(x,y)}{\Delta t} - \epsilon(u_{\theta,xx}+u_{\theta,yy}) + \beta_1 u_{\theta,x} + \beta_2 u_{\theta,y} + \sigma u_\theta - f(x,y,t_n).
		\]
		\For{each test function $v_{kl}(x,y)$}
		\State \textbf{Compute weak residual:}
		\[
		R_{kl}(\theta) := \int_\Omega R(x,y)\,v_{kl}(x,y)\,dxdy.
		\]
		\EndFor
		\State \textbf{Construct VPINN loss functional:}
		\[
		\mathcal{L}(\theta) := \frac{1}{M^2}\sum_{k,l=1}^M \big(R_{kl}(\theta)\big)^2.
		\]
		\State \textbf{Optimization:} Minimize $\mathcal{L}(\theta)$ using Adam, then fine-tune using L-BFGS.
		
		\State \textbf{Update:} $u^n(x,y) \gets u_\theta(x,y).$
		\EndFor
		\State \textbf{Output:} VPINN solution $\{u_\theta(x,y,t_n)\}_{n=1}^{N_t}$, and compute relative $L^2$ and $L^\infty$ error norms.
	\end{algorithmic}
\end{algorithm}

\begin{example}\label{example3}
	Consider the following one-parameter singularly perturbed IBVPs on G = $\Omega$$\times$$(0, T]:$
	\begin{equation}\label{E8}
		\left\{
		\begin{array} {ll}
			\frac{\partial u}{\partial t}-\epsilon \Delta u - u_x -u_y + u = f(x,y,t), & (x,y,t) \in G\\[8pt]
			u(x,y,0)=u_0(x,y), &  (x,y) \in \bar{\Omega},\\[8pt]
			u(x,y,t)=0, & (x,y,t) \in \partial \Omega\times(0,T].
		\end{array}\right.
	\end{equation}
\end{example}
The source term \( f(x,y,t) \) and initial condition \( u_0(x,y) \) 
are derived from the exact solution
\[
u(x,y,t) = e^{-t}\bigg(\sin(1-x)\sin(1-y)(1-e^{-x/\epsilon})(1-e^{-y/\epsilon})\bigg).\]

\begin{figure}[htbp]
	\centering
	\begin{subfigure}[b]{0.9\textwidth}
		\includegraphics[width=\textwidth]{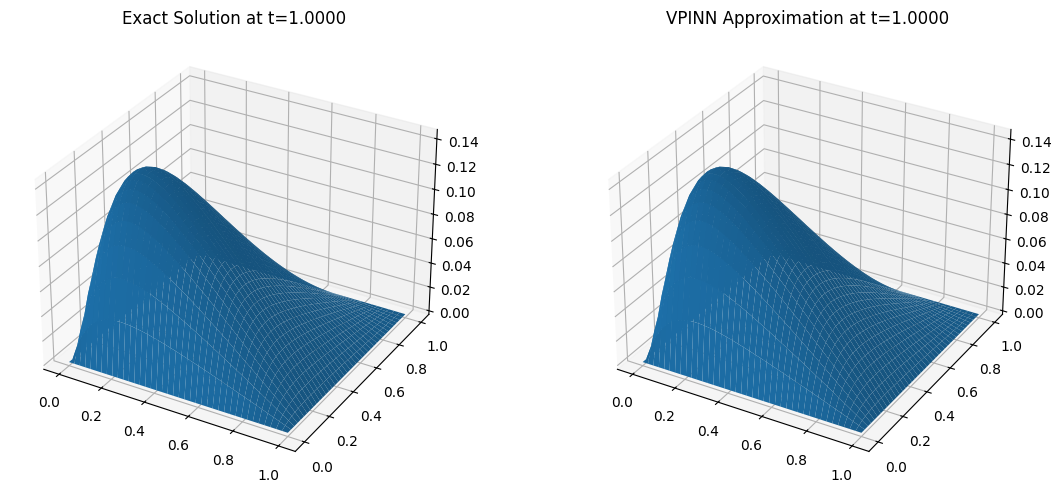}
		\caption{  $\epsilon=10^{-1}$}
		\label{fig:subfig7A}
	\end{subfigure}
	\hfill
	\begin{subfigure}[b]{0.9\textwidth}
		\includegraphics[width=\textwidth]{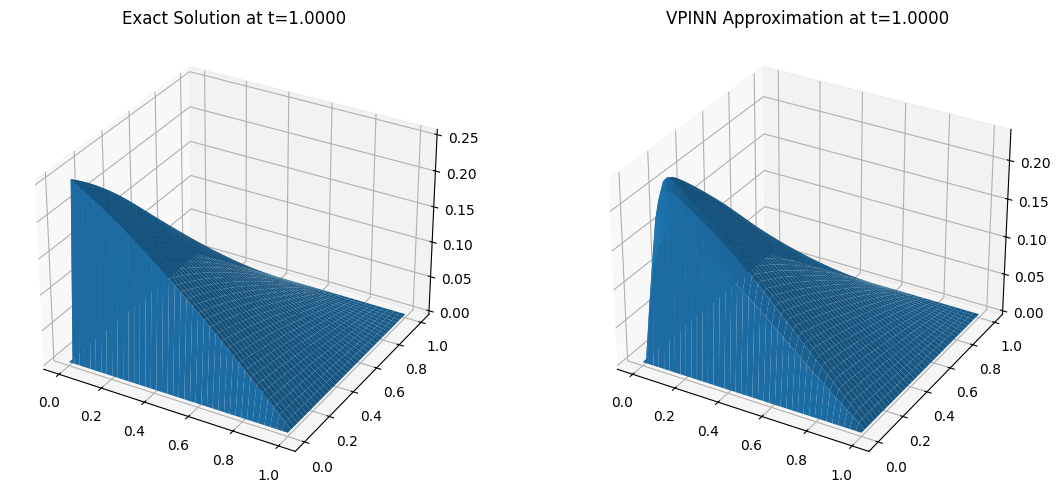}
		\caption{  $\epsilon=10^{-3}$}
		\label{fig:subfig7B}
	\end{subfigure}
	\hfill
	\begin{subfigure}[b]{0.9\textwidth}
		\includegraphics[width=\textwidth]{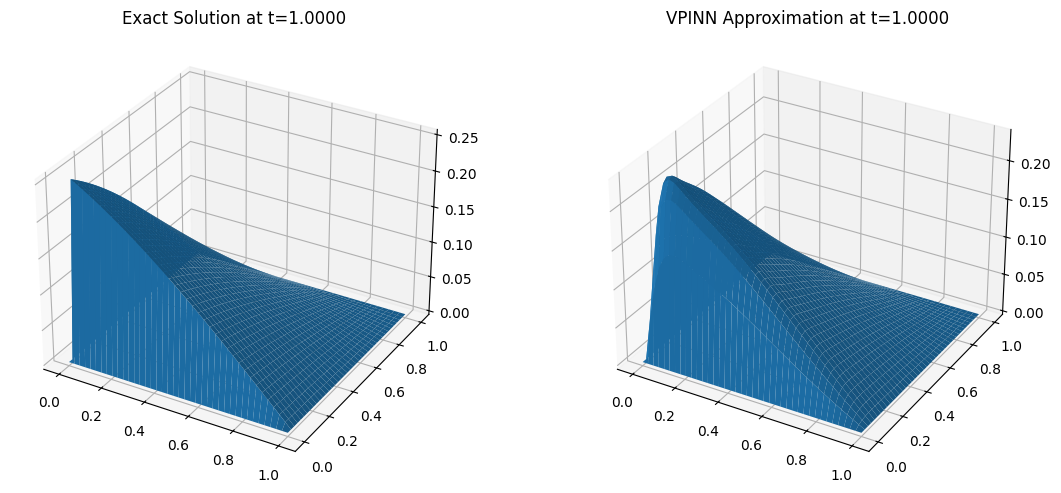}
		\caption{ $\epsilon=10^{-5}$}
		\label{fig:subfig7C}
	\end{subfigure}
	\caption{Comparison of the analytical and neural network solutions for test problem~\ref{example3} corresponding to various values of the perturbation parameters at the final time $T=1$.}\label{figure7}
\end{figure}

\begin{figure}[htbp]
	\centering
	
	\begin{subfigure}[b]{0.45\textwidth}
		\includegraphics[width=\textwidth,]{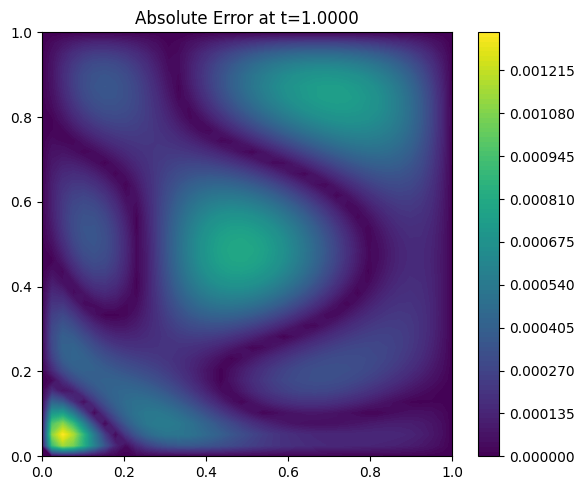}
		\caption{ $\epsilon=10^{-1}$}
		\label{fig:subfig8A}
	\end{subfigure}
	\hfill
	\begin{subfigure}[b]{0.45\textwidth}
		\includegraphics[width=\textwidth,]{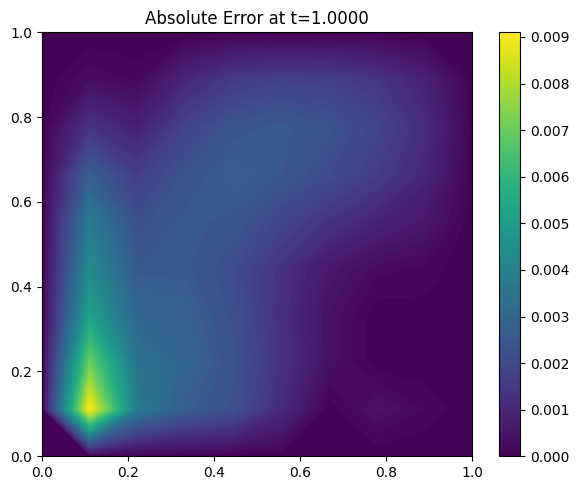}
		\caption{ $\epsilon=10^{-3}$}
		\label{fig:subfig8B}
	\end{subfigure}
	\hfill
	\begin{subfigure}[b]{0.45\textwidth}
		\includegraphics[width=\textwidth,]{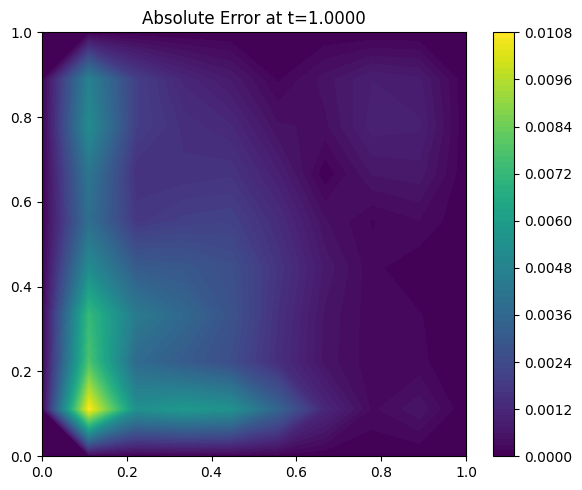}
		\caption{ $\epsilon=10^{-5}$}
		\label{fig:subfig8C}
	\end{subfigure}
	\caption{Absolute error at the final time $T=1$ for test problem~\ref{example3} corresponding to various values of the perturbation parameters.}\label{figure8}
\end{figure}

\begin{figure}[htbp]
	\centering
	\begin{subfigure}[b]{0.32\textwidth}
		\includegraphics[width=\textwidth,]{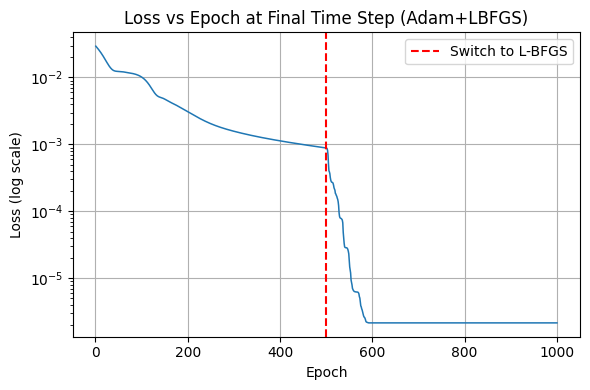}
		\caption{$\epsilon=10^{-1}$}
		\label{fig:subfig9A}
	\end{subfigure}
	\hfill
	\begin{subfigure}[b]{0.32\textwidth}
		\includegraphics[width=\textwidth,]{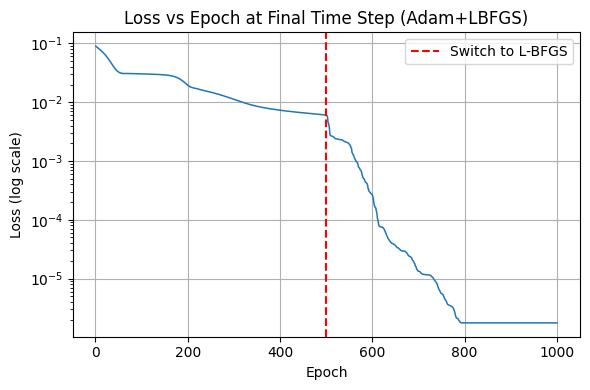}
		\caption{$\epsilon=10^{-3}$}
		\label{fig:subfig9B}
	\end{subfigure}
	\hfill
	\begin{subfigure}[b]{0.32\textwidth}
		\includegraphics[width=\textwidth,]{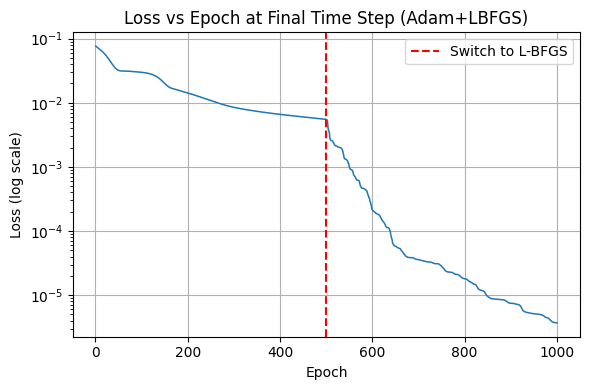}
		\caption{$\epsilon=10^{-5}$}
		\label{fig:subfig9C}
	\end{subfigure}
	\caption{Training loss versus epoch at the final time $T=1$ for test problem~\ref{example3} with corresponding to  various values of the perturbation parameters.}\label{figure9}
\end{figure}

The algorithm is executed for $1000$ iterations. The final errors at $T=1$ are summarized in Table~\ref{tbl_3}. Figure~\ref{figure7} compares the analytical and neural network solutions for different perturbation parameters at $T=1$, while the corresponding pointwise distributions of absolute errors are shown in Figure~\ref{figure8}. Figure~\ref{figure9} shows the loss versus epoch curves, illustrating the convergence behavior of the algorithm for various perturbation parameters at $T = 1$.

\begin{table}[htbp]
	\caption{\it{  Illustration of training loss and numerical error evolution for test problem~\ref{example3} at time $T=1$.}}\label{tbl_3}  
	\begin{tabular}{@{}lllllll@{}}
		\multicolumn 1 {c}{}  & \multicolumn 5 {c}
		{}\\
		\hline
		$\epsilon$ &Training Loss & Max Error & Rel-Max Error & $L_2$ Error & Rel-$L_2$ Error\\
		\hline
		$10^{-1}$  & 2.1600e-06 & 1.3279e-03 &  9.2001e-03 & 3.2441e-04 & 4.9370e-03\\ [6pt]
		
		$10^{-2}$  & 5.7879e-06 & 8.9176e-03 &  4.0217e-02 & 2.1246e-03 & 2.7392e-02 \\ [6pt]
		
		$10^{-3}$  & 1.7581e-06 &  9.0439e-03 & 4.0786e-02 &  1.9694e-03 & 2.5390e-02 \\ [6pt]
		
		$10^{-4}$    &  1.8143e-05 &  1.6519e-02  & 7.4496e-02 &  2.9020e-03 &   3.7413e-02  \\ [6pt]
		
		$10^{-5}$    & 3.6927e-06 & 1.0685e-02  & 4.8188e-02 &   2.4062e-03 & 3.1021e-02 \\ [6pt]
		
		\hline
	\end{tabular}
\end{table}

\subsection{Two-dimensional two-parameter elliptic BVPs}\label{sec3.3}
Consider the two-dimensional two-parameter elliptic BVPs:
\begin{equation}\label{E9}
	\left\{
	\begin{array} {ll}
		-\epsilon \Delta u + \mu\boldsymbol{\beta}(x,y)\cdot\nabla u + \sigma(x,y) u = f, & (x,y) \in \Omega=(0,\,1)\times(0,\,1), \\[8pt]
		u = 0, &  (x,y) \in \partial\Omega.\\[8pt]
	\end{array}\right.
\end{equation}
Here, $0<\epsilon, \mu \ll 1$ are small parameters, $\boldsymbol{\beta}(x,y) = (\beta_1, \beta_2)^{T}$ with $\beta_1 \geq  \alpha_1 > 0$, $\beta_2 \geq  \alpha_2 > 0$,
and $\sigma(x,y) \geq \delta > 0$ for any $(x,y)\in\overline{\Omega}$. 
Furthermore, we assume the stability condition:
\[
c_0^2 = \left(\sigma - \tfrac{1}{2}\mu\nabla\cdot \boldsymbol{\beta}\right) \geq \omega > 0.
\]
The functions $\boldsymbol{\beta}$, $\sigma$, and $f$ are assumed to be sufficiently smooth, and $f$ satisfies the following compatibility condition:
\[f(0,0)=f(0,1)=f(1,0)=f(1,1)=0.\]

The finite element discretization and corresponding VPINN algorithm can be constructed analogously to the model problem defined in equation~\eqref{E1}.

Two test problems are considered: one involving constant coefficients and another involving variable coefficients.
\begin{example}\label{example4}
	Consider the $2$D two-parameter elliptic BVP:
	\begin{equation}\label{E10}
		\left\{
		\begin{array} {ll}
			-\epsilon \Delta u + \mu u_x +  u = f, & (x,y) \in \Omega= (0, 1)\times (0,1), \\[8pt]
			u = 0, &  (x,y) \in \partial\Omega.\\[8pt]
		\end{array}\right.
	\end{equation}
	The source function $f$ is derived from the exact solution
	\begin{eqnarray}
		u(x,y) = &0.25&\times 
		\Big(1 - \exp\big(-\tfrac{\mu l_1 x}{2 \epsilon}\big)\Big)\times 
		\Big(1 - \exp\big(-\tfrac{y}{\sqrt{\epsilon}}\big)\Big)\nonumber \\
		&&\times\Big(1 - \exp\big(-\tfrac{\mu l_2 (1-x)}{2 \epsilon}\big)\Big) \times
		\Big(1 - \exp\big(-\tfrac{1-y}{\sqrt{\epsilon}}\big)\Big),\nonumber
	\end{eqnarray}
	
	where $l_{1,2} = \bigg(\sqrt{\big(1+ 16 \frac{\epsilon}{\mu^2}\big)}\mp 1\bigg)$.
	
\end{example}

\begin{example}\label{example5}
	Consider the $2$D two-parameter elliptic BVP:
	\begin{equation}\label{E11}
		\left\{
		\begin{array} {ll}
			-\epsilon \Delta u +\mu(3-x) u_x  +  u = f, & (x,y) \in \Omega= (0, 1)\times (0,1), \\[8pt]
			u = 0, &  (x,y) \in \partial\Omega.\\[8pt]
		\end{array}\right.
	\end{equation}
	The source function $f$ is derived from the exact solution
	\begin{eqnarray}
		u(x,y) = \bigg(1+\frac{\sin(8x)}{2}\bigg)\times v(x,y),\nonumber 
	\end{eqnarray}
	here, 
	\begin{eqnarray}
		v(x,y) = &0.25&\times 
		\Big(1 - \exp\big(-\tfrac{\mu l_1 x}{2 \epsilon}\big)\Big)\times 
		\Big(1 - \exp\big(-\tfrac{y}{\sqrt{\epsilon}}\big)\Big)\nonumber \\
		&&\times\Big(1 - \exp\big(-\tfrac{\mu l_2 (1-x)}{2 \epsilon}\big)\Big) \times
		\Big(1 - \exp\big(-\tfrac{1-y}{\sqrt{\epsilon}}\big)\Big).\nonumber
	\end{eqnarray}
	and $l_{1,2}=\bigg(\sqrt{\big(1+ 16 \frac{\epsilon}{\mu^2}\big)}\mp 1\bigg)$. 
\end{example}

The algorithm is executed for $1000$ iterations. The resulting maximum error, $L_2$ error, and loss at the final iteration are reported in Table~\ref{tbl_4} for test problem~\ref{example4} and Table~\ref{tbl_5} for test problem~\ref{example5}.  

\begin{figure}[htbp]
	\centering
	\begin{subfigure}[b]{0.9\textwidth}
		\includegraphics[width=\textwidth]{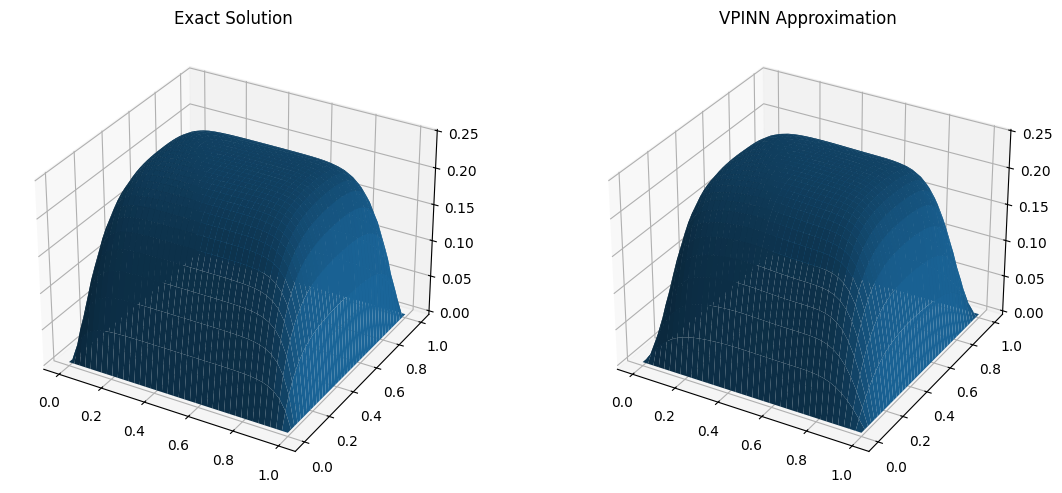}
		\caption{ $\epsilon = 10^{-2}$ and $\mu = 10^{-3}$}
		\label{fig:subfig10A}
	\end{subfigure}
	\hfill
	\begin{subfigure}[b]{0.9\textwidth}
		\includegraphics[width=\textwidth]{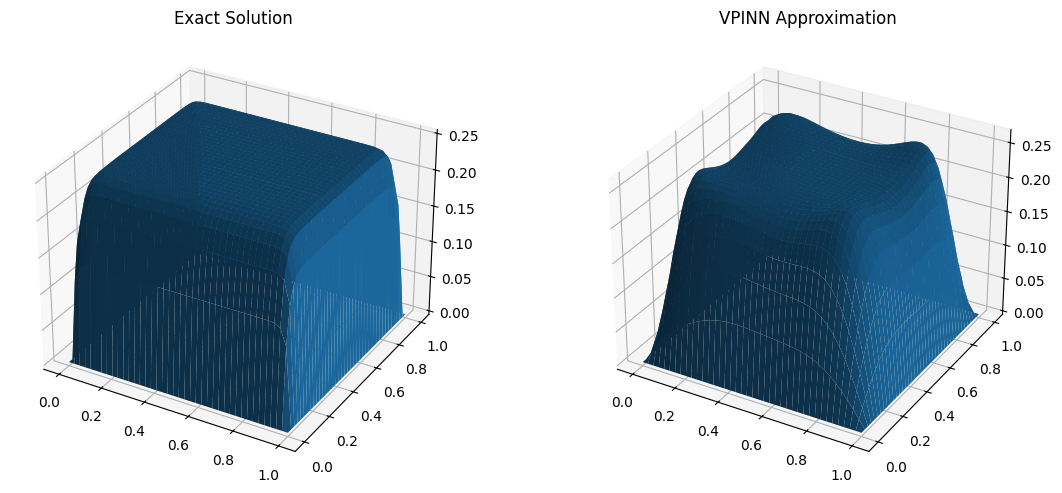}
		\caption{ $\epsilon = 10^{-3}$ and $\mu = 10^{-4}$}
		\label{fig:subfig10B}
	\end{subfigure}
	\hfill
	\begin{subfigure}[b]{0.9\textwidth}
		\includegraphics[width=\textwidth]{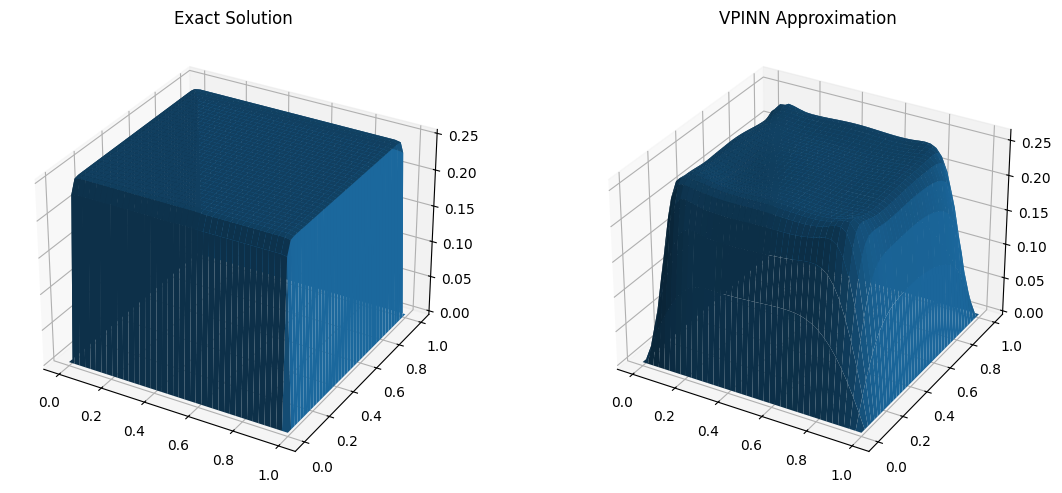}
		\caption{ $\epsilon = 10^{-4}$ and $\mu = 10^{-5}$}
		\label{fig:subfig10C}
	\end{subfigure}
	\caption{Comparison of the analytical and neural network solutions for test problem~\ref{example4} corresponding to various values of the perturbation parameters.}\label{figure10}
\end{figure}

The comparison between the  analytical and neural network solutions for different parameter values is shown in Figure~\ref{figure10} for test problem~\ref{example4} and Figure~\ref{figure13} for test problem~\ref{example5}.  

\begin{figure}[htbp]
	\centering
	
	\begin{subfigure}[b]{0.45\textwidth}
		\includegraphics[width=\textwidth,]{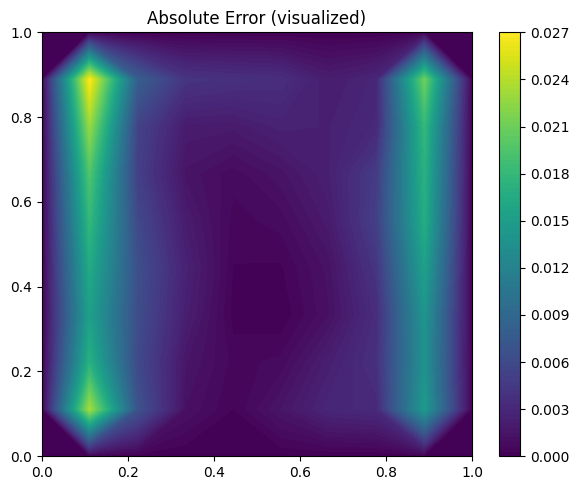}
		\caption{ $\epsilon = 10^{-2}$ and $\mu = 10^{-3}$}
		\label{fig:subfig11A}
	\end{subfigure}
	\hfill
	\begin{subfigure}[b]{0.45\textwidth}
		\includegraphics[width=\textwidth,]{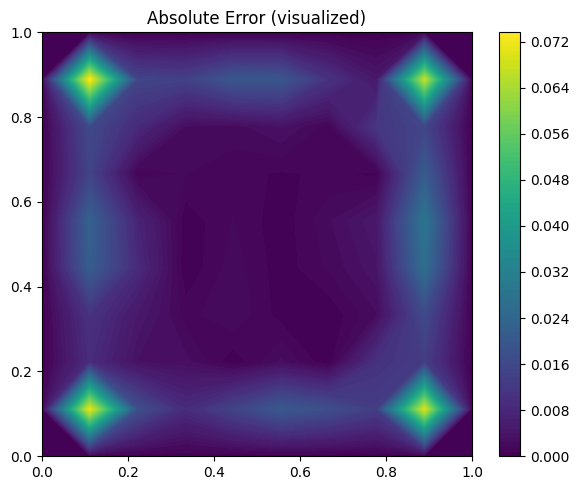}
		\caption{ $\epsilon = 10^{-3}$ and $\mu = 10^{-4}$}
		\label{fig:subfig11B}
	\end{subfigure}
	\hfill
	\begin{subfigure}[b]{0.45\textwidth}
		\includegraphics[width=\textwidth,]{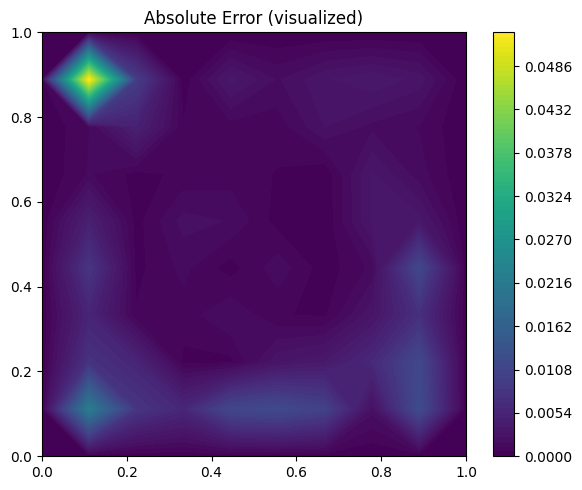}
		\caption{ $\epsilon = 10^{-4}$ and $\mu = 10^{-5}$}
		\label{fig:subfig11C}
	\end{subfigure}
	\caption{Absolute error for test problem~\ref{example4} corresponding to various values of the perturbation parameters.}\label{figure11}
\end{figure}

\begin{figure}[htbp]
	\centering
	
	\begin{subfigure}[b]{0.32\textwidth}
		\includegraphics[width=\textwidth,]{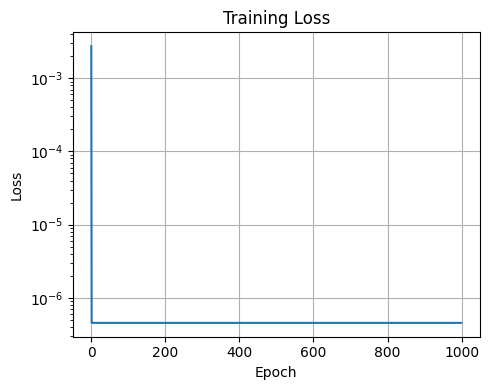}
		\caption{ $\epsilon = 10^{-2}$ and $\mu = 10^{-3}$}
		\label{fig:subfig12A}
	\end{subfigure}
	\hfill
	\begin{subfigure}[b]{0.32\textwidth}
		\includegraphics[width=\textwidth,]{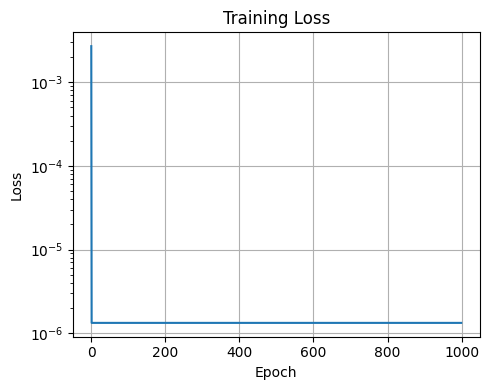}
		\caption{ $\epsilon = 10^{-3}$ and $\mu = 10^{-4}$}
		\label{fig:subfig12B}
	\end{subfigure}
	\hfill
	\begin{subfigure}[b]{0.32\textwidth}
		\includegraphics[width=\textwidth,]{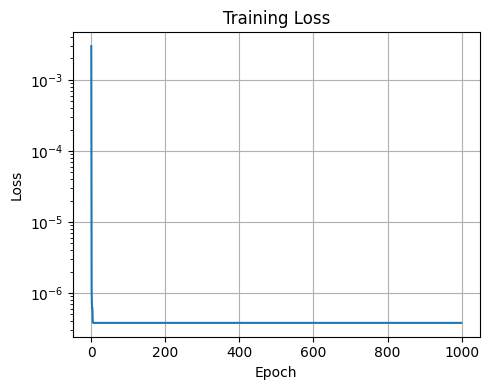}
		\caption{ $\epsilon = 10^{-4}$ and $\mu = 10^{-5}$}
		\label{fig:subfig12C}
	\end{subfigure}
	\caption{Training loss versus epoch for test problem~\ref{example4} corresponding to various values of the perturbation parameters.}\label{figure12}
\end{figure}

\begin{figure}[htbp]
	\centering
	\begin{subfigure}[b]{0.9\textwidth}
		\includegraphics[width=\textwidth]{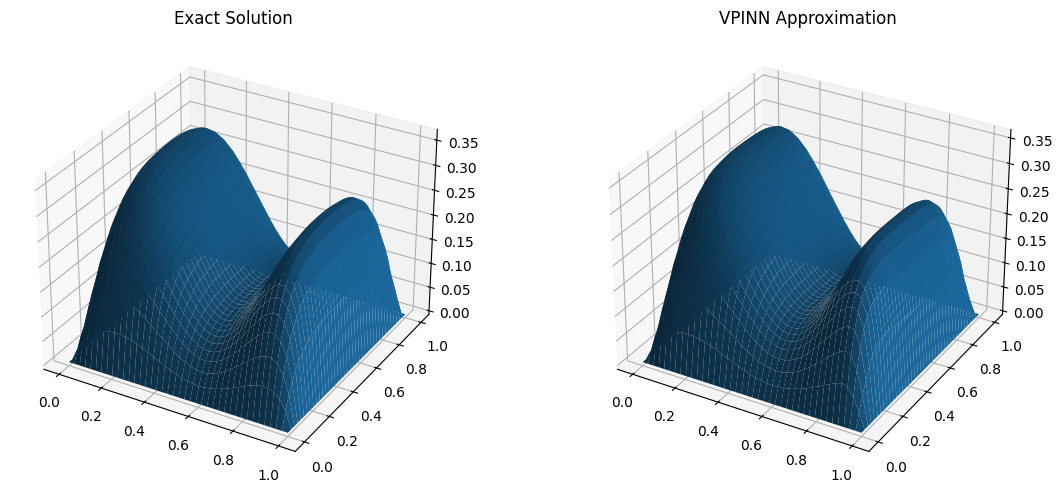}
		\caption{ $\epsilon = 10^{-2}$ and $\mu = 10^{-3}$}
		\label{fig:subfig13A}
	\end{subfigure}
	\hfill
	\begin{subfigure}[b]{0.9\textwidth}
		\includegraphics[width=\textwidth]{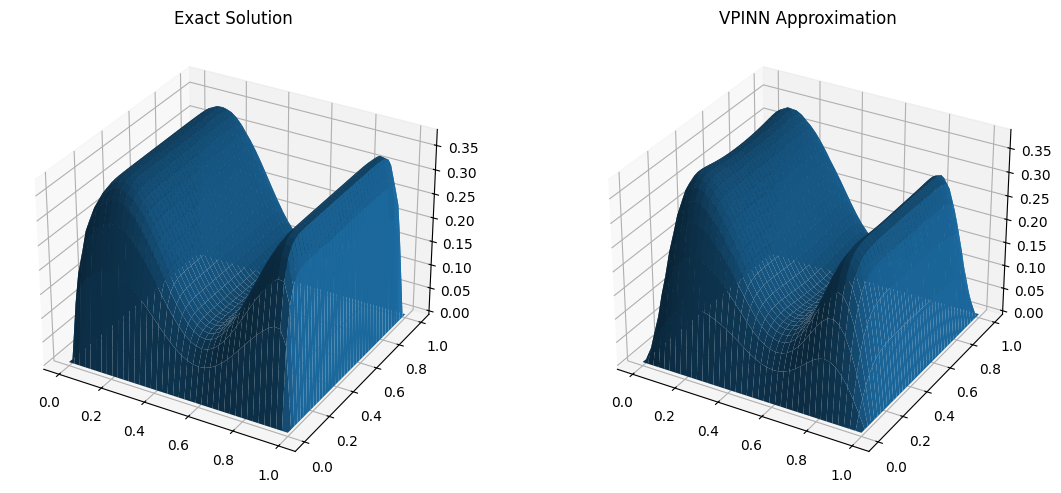}
		\caption{ $\epsilon = 10^{-3}$ and $\mu = 10^{-4}$}
		\label{fig:subfig13B}
	\end{subfigure}
	\hfill
	\begin{subfigure}[b]{0.9\textwidth}
		\includegraphics[width=\textwidth]{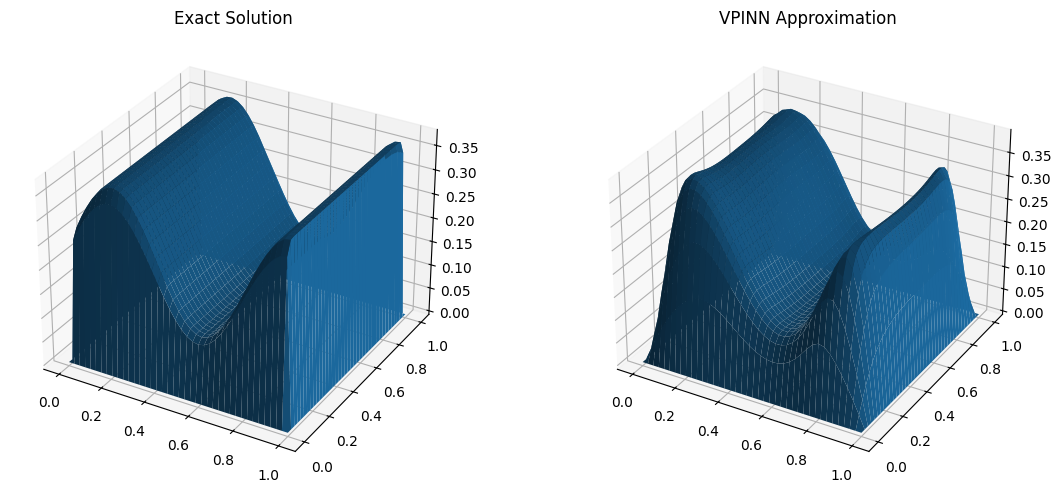}
		\caption{ $\epsilon = 10^{-4}$ and $\mu = 10^{-5}$}
		\label{fig:subfig13C}
	\end{subfigure}
	\caption{Comparison of the analytical and neural network solutions for test problem~\ref{example5} corresponding to various values of the perturbation parameters.}\label{figure13}
\end{figure}

\begin{figure}[htbp]
	\centering
	
	\begin{subfigure}[b]{0.45\textwidth}
		\includegraphics[width=\textwidth,]{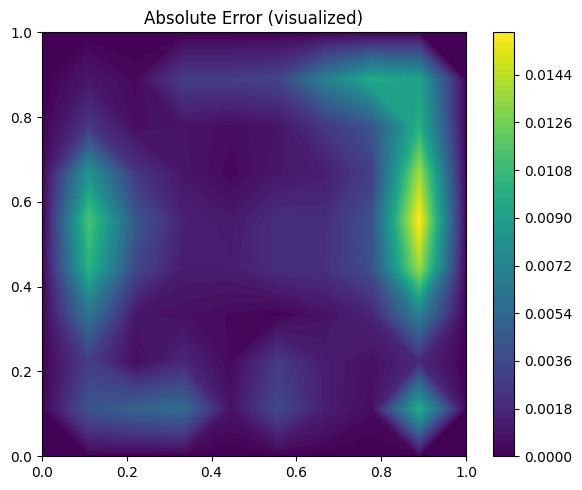}
		\caption{ $\epsilon = 10^{-2}$ and $\mu = 10^{-3}$}
		\label{fig:subfig14A}
	\end{subfigure}
	\hfill
	\begin{subfigure}[b]{0.45\textwidth}
		\includegraphics[width=\textwidth,]{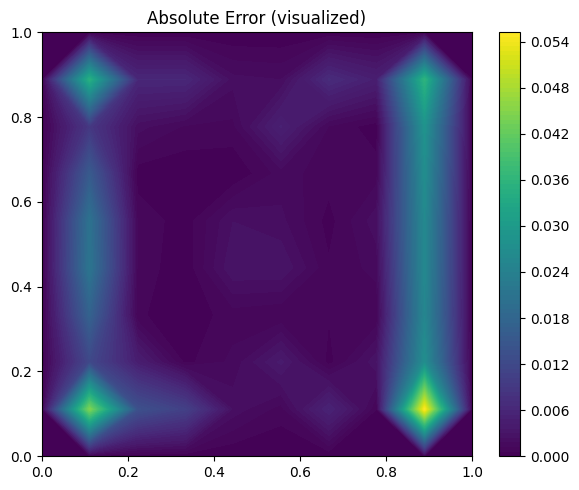}
		\caption{ $\epsilon = 10^{-3}$ and $\mu = 10^{-4}$}
		\label{fig:subfig14B}
	\end{subfigure}
	\hfill
	\begin{subfigure}[b]{0.45\textwidth}
		\includegraphics[width=\textwidth,]{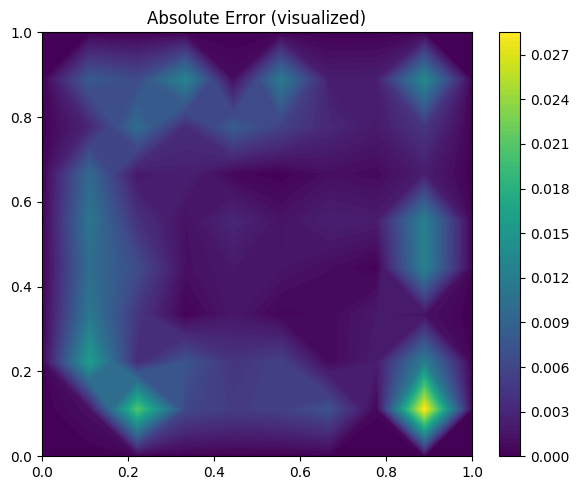}
		\caption{ $\epsilon = 10^{-4}$ and $\mu = 10^{-5}$}
		\label{fig:subfig14C}
	\end{subfigure}
	\caption{Absolute error for test problem~\ref{example5} corresponding to various values of the perturbation parameters.}\label{figure14}
\end{figure}

\begin{figure}[htbp]
	\centering
	
	\begin{subfigure}[b]{0.32\textwidth}
		\includegraphics[width=\textwidth,]{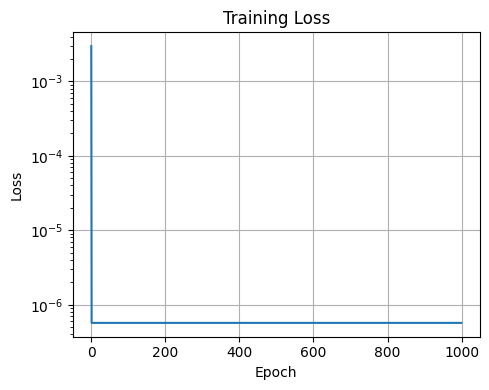}
		\caption{ $\epsilon = 10^{-2}$ and $\mu = 10^{-3}$}
		\label{fig:subfig15A}
	\end{subfigure}
	\hfill
	\begin{subfigure}[b]{0.32\textwidth}
		\includegraphics[width=\textwidth,]{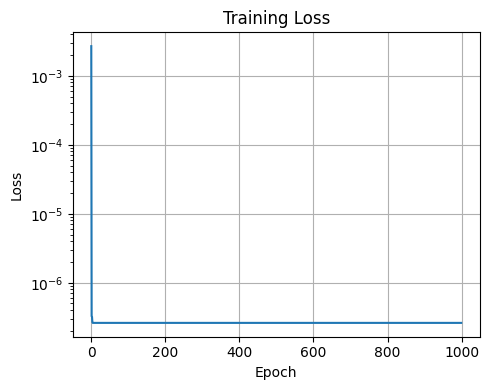}
		\caption{ $\epsilon = 10^{-3}$ and $\mu = 10^{-4}$}
		\label{fig:subfig15B}
	\end{subfigure}
	\hfill
	\begin{subfigure}[b]{0.32\textwidth}
		\includegraphics[width=\textwidth,]{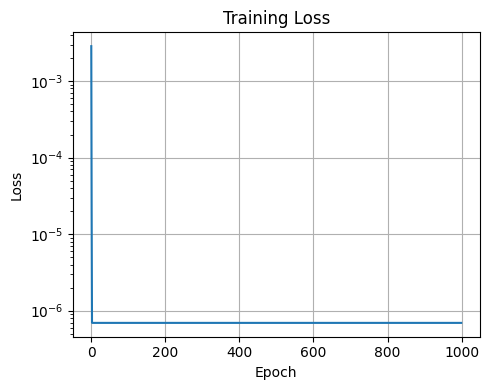}
		\caption{ $\epsilon = 10^{-4}$ and $\mu = 10^{-5}$}
		\label{fig:subfig15C}
	\end{subfigure}
	\caption{Training loss versus epoch for test problem~\ref{example5} corresponding to various values of  the perturbation parameters.}\label{figure15}
\end{figure}

The absolute error plots are presented in Figure~\ref{figure11} for test problem~\ref{example4} and Figure~\ref{figure14} for test problem~\ref{example5}, while the training loss versus epoch curves for the final iteration are shown in Figure~\ref{figure12} and  Figure~\ref{figure15} corresponding to test problem~\ref{example4} and test problem~\ref{example5}.

\begin{table}[htbp]
	\caption{\it{ Illustration of training loss and numerical error evolution for test problem~\ref{example4}.}}\label{tbl_4}  
	\begin{tabular}{@{}lllllll@{}}
		\multicolumn 1 {c}{}  & \multicolumn 5 {c}
		{}\\
		\hline
		$\epsilon,\ \ \quad \mu$ &Training Loss & Max Error & Rel-Max Error & $L_2$ Error & Rel-$L_2$ Error \\
		\hline
		$10^{-1},10^{-2}$  &  2.5401e-07 & 9.6440e-04 & 6.7707e-03 & 2.5607e-04 &  3.2958e-03 \\ [6pt]
		$10^{-2},10^{-3}$  &  4.5312e-07 &  2.6850e-02 &  8.9031e-02 &  7.7171e-03 & 4.4860e-02 \\ [6pt]

		$10^{-3},10^{-4}$ &   1.3252e-06 & 7.3284e-02  & 8.2121e-02 &   1.6425e-02 & 8.2773e-02 \\ [6pt]
		
		$10^{-4},10^{-5}$ & 3.8090e-07  & 5.2750e-02  &  8.3948e-02 &   6.8092e-03 &   3.4046e-02 \\ [6pt]
		\hline
	\end{tabular}
\end{table}

\begin{table}[htbp]
	\caption{\it{  Illustration of training loss and numerical error evolution for test problem~\ref{example5}.}}\label{tbl_5}  
	\begin{tabular}{@{}lllllll@{}}
		\multicolumn 1 {c}{}  & \multicolumn 5 {c}
		{}\\
		\hline
		$\epsilon,\ \ \quad \mu$ &Training Loss & Max Error & Rel-Max Error & $L_2$ Error & Rel-$L_2$ Error \\
		\hline
		$10^{-1},10^{-2}$  &  1.1794e-06 & 2.5570e-03  & 1.4727e-02 &  7.0438e-04 & 8.9877e-03\\ [6pt]
		$10^{-2},10^{-3}$  &  5.7021e-07 &  1.5983e-02 &  4.4130e-02 & 4.1249e-03 & 2.2107e-02 \\ [6pt]
		
		$10^{-3},10^{-4}$ &    2.6163e-07 & 5.5043e-02  &  5.3746e-02 &  1.1758e-02 &  5.3746e-02  \\ [6pt]
		
		$10^{-4},10^{-5}$ &  6.9310e-07 & 2.8340e-02  &  7.6114e-02 &  5.9708e-03 &   2.7077e-02 \\ [6pt]
		\hline
	\end{tabular}
\end{table}

\subsection{Two-dimensional parabolic two-parameter IBVPs }\label{sec3.4}
Consider the  two-dimensional two-parameters singularly perturbed parabolic IBVPs defined on the domain $G$ = $\Omega$ $\times$ $(0, T]$, where $\Omega = (0, 1)\times (0,1):$
\begin{equation}\label{E12}
	\left\{
	\begin{array} {ll}
		\frac{\partial u}{\partial t}-\epsilon \Delta u + \mu\boldsymbol{\beta}(x,y)\cdot\nabla u + \sigma(x,y) u = f(x,y,t), & (x,y,t) \in G, \\[8pt]
		u(x,y,0)=u_0(x,y), &  (x,y) \in \bar{\Omega},\\[8pt]
		u(x,y,t)=0 & (x,y,t) \in \partial \Omega\times(0,T].
	\end{array}\right.
\end{equation}

Here, $0<\epsilon\ll 1$ and $0<\mu\ll 1$ are small parameters, $ \boldsymbol{\beta}, \sigma$ and $f$ are smooth and $\boldsymbol{\beta}(x,y) = (\beta_1, \beta_2)^{T}$ with $\beta_1 \geq  \alpha_1 > 0$, $\beta_2 \geq  \alpha_2 > 0$,
and $\sigma(x,y) \geq \delta > 0$ for any $(x,y)\in\overline{\Omega}$. 
Furthermore, we assume the stability condition
\[
c_0^2 := \left(\sigma - \tfrac{1}{2}\mu \nabla\cdot \boldsymbol{\beta}\right) \geq \omega > 0.
\]

The corresponding finite element formulation and numerical algorithm follow analogously from the model problem~(\ref{E5}).

\begin{example}\label{example6}
	Consider the following two-parameter singularly perturbed IBVPs on G = $\Omega$$\times$$(0, T]:$
	\begin{equation}\label{E13}
		\left\{
		\begin{array} {ll}
			\frac{\partial u}{\partial t}-\epsilon \Delta u  + \mu u_x +  u  = f(x,y,t), & (x,y,t) \in G, \\[8pt]
			u(x,y,0)=u_0(x,y), &  (x,y) \in \bar{\Omega},\\[8pt]
			u(x,y,t)=0, & (x,y,t) \in \partial \Omega\times(0,T].
		\end{array}\right.
	\end{equation}
\end{example}
The initial condition \( u_0(x,y) \) and source term \( f(x,y,t) \) 
are chosen from the exact solution
\[
u(x,y,t) = e^{-t}\bigg(\bigg(1+\frac{\sin(8x)}{2}\bigg)\times v(x,y)\bigg),\]

where 
\begin{eqnarray}
	v(x,y) = &0.25&\times 
	\Big(1 - \exp\big(-\tfrac{\mu l_1 x}{2 \epsilon}\big)\Big)\times 
	\Big(1 - \exp\big(-\tfrac{y}{\sqrt{\epsilon}}\big)\Big)\nonumber \\
	&&\times\Big(1 - \exp\big(-\tfrac{\mu l_2 (1-x)}{2 \epsilon}\big)\Big) \times
	\Big(1 - \exp\big(-\tfrac{1-y}{\sqrt{\epsilon}}\big)\Big).\nonumber
\end{eqnarray}
and $l_{1,2}=\bigg(\sqrt{\big(1+ 16 \frac{\epsilon}{\mu^2}\big)}\mp 1\bigg)$.

\begin{figure}[htbp]
	\centering
	\begin{subfigure}[b]{0.9\textwidth}
		\includegraphics[width=\textwidth]{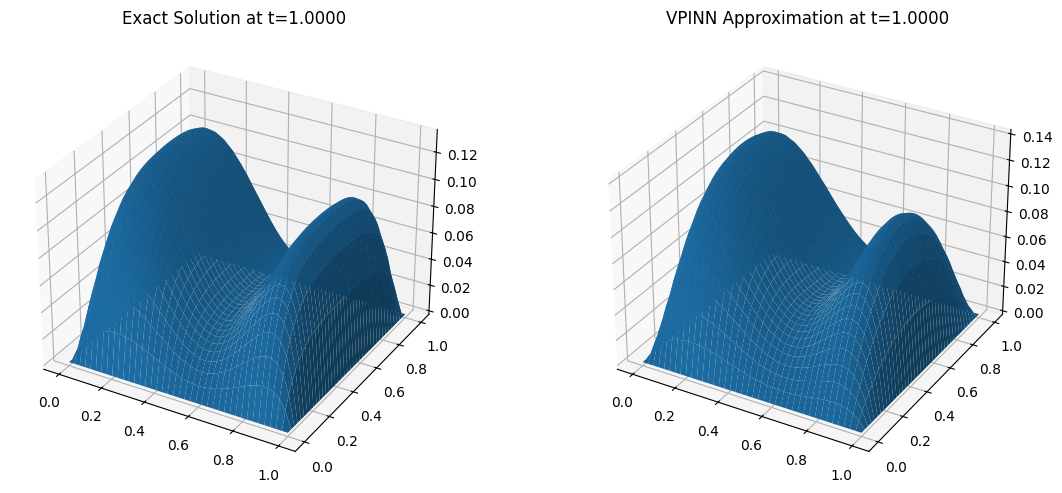}
		\caption{ $\epsilon = 10^{-2}$ and $\mu = 10^{-3}$}
		\label{fig:subfig16A}
	\end{subfigure}
	\hfill
	\begin{subfigure}[b]{0.9\textwidth}
		\includegraphics[width=\textwidth]{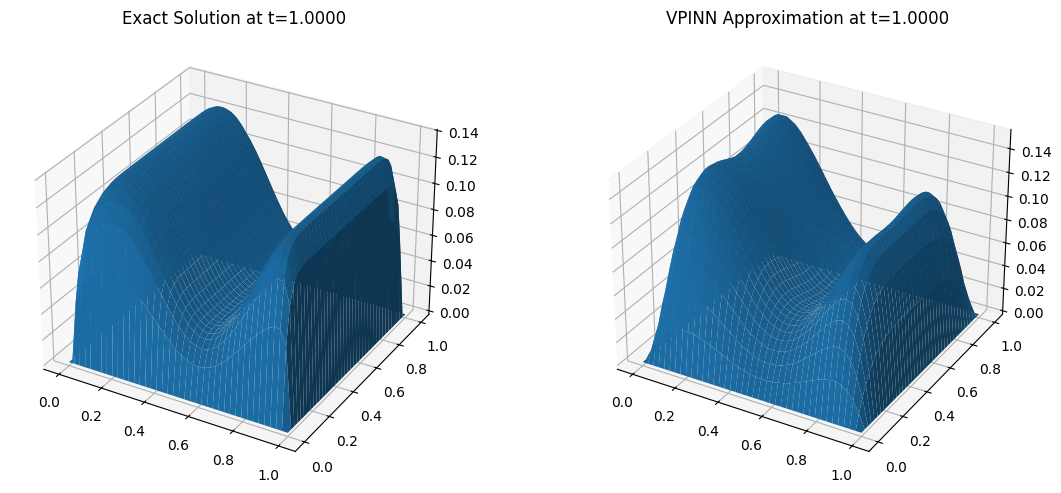}
		\caption{ $\epsilon = 10^{-3}$ and $\mu = 10^{-4}$}
		\label{fig:subfig16B}
	\end{subfigure}
	\hfill
	\begin{subfigure}[b]{0.9\textwidth}
		\includegraphics[width=\textwidth]{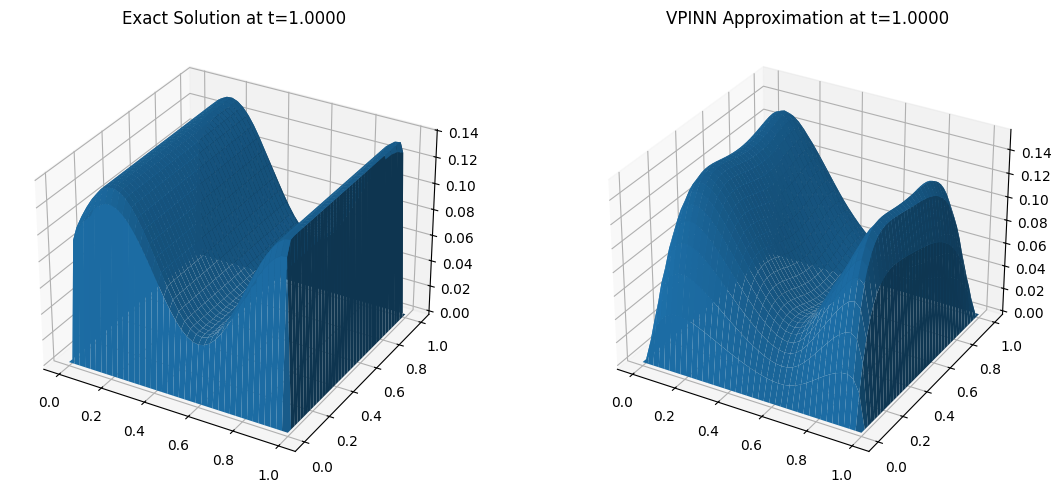}
		\caption{ $\epsilon = 10^{-4}$ and $\mu = 10^{-5}$}
		\label{fig:subfig16C}
	\end{subfigure}
	\caption{Comparison of the analytical and neural network solutions for test problem~\ref{example6} corresponding to various values of the perturbation parameters at the final time $T=1$.}\label{figure16}
\end{figure}

\begin{figure}[htbp]
	\centering
	
	\begin{subfigure}[b]{0.45\textwidth}
		\includegraphics[width=\textwidth]{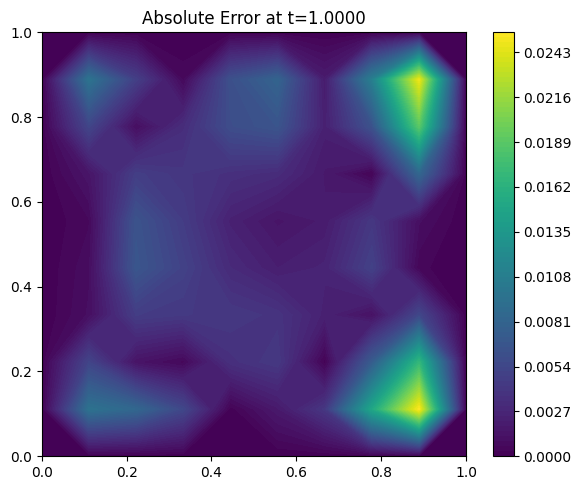}
		\caption{ $\epsilon = 10^{-2}$ and $\mu = 10^{-3}$}
		\label{fig:subfig17A}
	\end{subfigure}
	\hfill
	\begin{subfigure}[b]{0.45\textwidth}
		\includegraphics[width=\textwidth,]{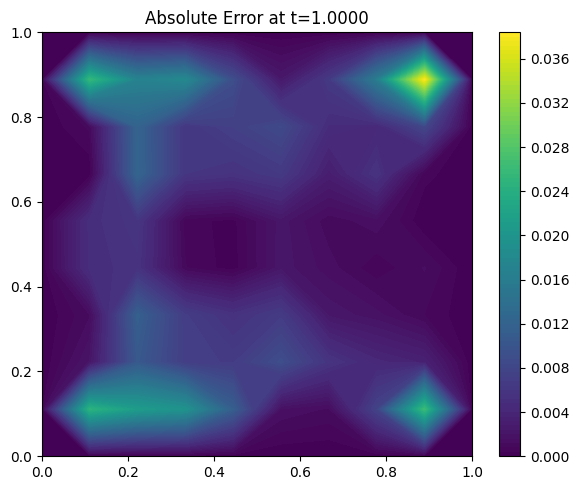}
		\caption{ $\epsilon = 10^{-3}$ and $\mu = 10^{-4}$}
		\label{fig:subfig17B}
	\end{subfigure}
	\hfill
	\begin{subfigure}[b]{0.45\textwidth}
		\includegraphics[width=\textwidth,]{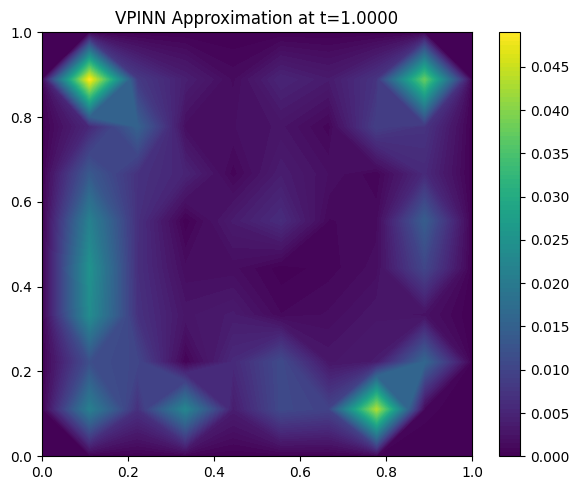}
		\caption{ $\epsilon = 10^{-4}$ and $\mu = 10^{-5}$}
		\label{fig:subfig17C}
	\end{subfigure}
	\caption{Absolute error for test problem~\ref{example6} corresponding to various values of  the perturbation parameters at the final time $T=1$.}\label{figure17}
\end{figure}

\begin{figure}[htbp]
	\centering
	\begin{subfigure}[b]{0.32\textwidth}
		\includegraphics[width=\textwidth]{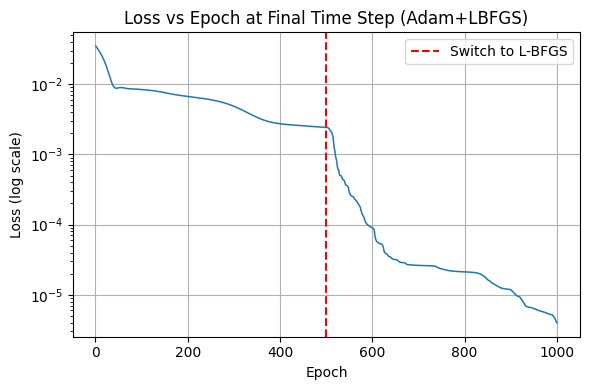}
		\caption{ $\epsilon = 10^{-2}$ and $\mu = 10^{-3}$}
		\label{fig:subfig18A}
	\end{subfigure}\hfill
	\begin{subfigure}[b]{0.32\textwidth}
		\includegraphics[width=\textwidth]{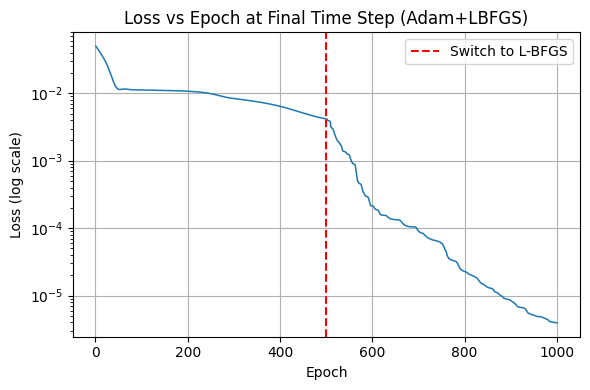}
		\caption{ $\epsilon = 10^{-3}$ and $\mu = 10^{-4}$}
		\label{fig:subfig18B}
	\end{subfigure}\hfill
	\begin{subfigure}[b]{0.32\textwidth}
		\includegraphics[width=\textwidth]{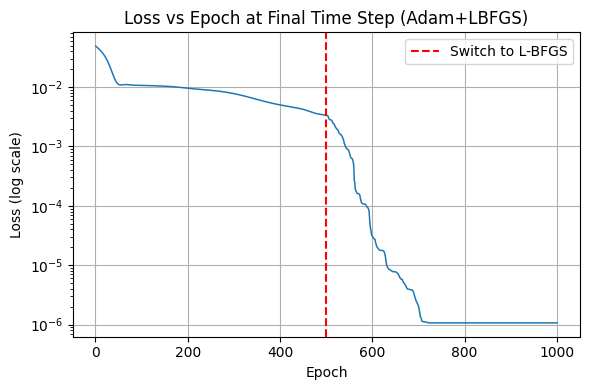}
		\caption{ $\epsilon = 10^{-4}$ and $\mu = 10^{-5}$}
		\label{fig:subfig18C}
	\end{subfigure}
	\caption{Training loss versus epoch for test problem~\ref{example6} corresponding to various values of the perturbation parameters at the final time $T=1$.}
	\label{figure18}
\end{figure}

The algorithm is executed at $1000$ iterations. The final errors for $T=1$ are summarized in Table~\ref{tbl_6}. Figure~\ref{figure16} compares the analytical and neural network solutions for various perturbation parameters at $T=1$, while the corresponding pointwise distributions of absolute errors are depicted in Figure~\ref{figure17}. Figure~\ref{figure18} presents the loss versus epoch curves, highlighting the convergence of the algorithm at the last iteration for various perturbation parameters at $T=1$.

\begin{table}[htbp]
	\caption{\it{ Illustration of training loss and numerical error evolution for test problem~\ref{example6} at time $T=1$.}}\label{tbl_6}  
	\begin{tabular}{@{}lllllll@{}}
		\multicolumn 1 {c}{}  & \multicolumn 5 {c}
		{}\\
		\hline
		$\epsilon, \quad \ \mu$ &Training Loss & Max Error & Rel-Max Error & $L_2$ Error & Rel-$L_2$ Error \\
		\hline
		$10^{-1},10^{-2}$  &  3.6371e-07 & 1.4507e-03  & 2.2713e-02 & 6.9817e-04 & 2.4216e-02 \\ [6pt]
		
		$10^{-2},10^{-3}$  &  3.9734e-06 &2.5320e-02 & 2.9003e-02 &  5.8064e-03 & 8.4592e-02 \\ [6pt]
		
		$10^{-3},10^{-4}$ &    3.9383e-06 & 3.8213e-02 &  4.8639e-02 &  7.1570e-03 &  5.2887e-02 \\ [6pt]
		
		$10^{-4},10^{-5}$ &  1.0703e-06 &  4.8758e-02  &   6.2020e-02 &    4.6401e-03 &   2.1042e-02 \\ [6pt]
		\hline
	\end{tabular}
\end{table}

\section{Conclusion}\label{Conclusion}
We present a VPINN methodology employing the Petrov-Galerkin framework for two-dimensional SPPs, with the trial solution space represented by DNNs, while tensor-product hat functions are employed as test functions in the weak formulation. This construction effectively captured sharp boundary layers in two spatial dimensions without the need for fine mesh refinement. The backward Euler method is employed for temporal discretization in time-dependent cases, ensuring stability in the stiff regime. The training strategy combined the efficiency of Adam optimizer during the early stages with the robustness of L-BFGS for fine-tuning.  Dirichlet boundary conditions are imposed directly, while the source terms are computed by automatic differentiation, maintaining high accuracy for the variational framework. The proposed two-dimensional VPINN method exhibited stable and accurate performance, providing a mesh-free and flexible solution approach for steady-state and time-dependent SPPs.

\section*{Funding}
The first author gratefully acknowledges the financial support provided by the National Institute of Technology, Tiruchirappalli, for his research.

\section*{Conflict of Interest}
The authors declare that they have no conflicts of interest.

\section*{Data Availability}
Not applicable.

\end{document}